\documentclass{cicp}
\usepackage{tabularx}
\usepackage{makecell}
\usepackage{mathrsfs}
\usepackage{amsmath}
\usepackage{stmaryrd}
\usepackage{bbding}
\usepackage{dcolumn}
\usepackage{graphicx}
\usepackage{amsfonts}
\usepackage{amssymb}
\usepackage{psfrag}
\usepackage{wrapfig}
\usepackage{subfigure}
\usepackage{makeidx}
\usepackage{bm}
\usepackage{epsf}
\usepackage{color}
\usepackage{epsfig}
\usepackage{setspace}
\usepackage{graphicx}
\usepackage{epstopdf}
\usepackage{psfrag}
\usepackage{subfigure}
\usepackage{multirow}
\usepackage{diagbox}
\usepackage{verbatim}

\usepackage{makecell}
\usepackage{float}
\usepackage{esint}
\usepackage{comment}
\usepackage{movie15}
\usepackage{hyperref}
\hypersetup{hidelinks}
\usepackage[ruled,linesnumbered]{algorithm2e}
\allowdisplaybreaks[4]

\begin{document}
\title{An Adaptive Reconstruction Method for Arbitrary High-Order Accuracy Using Discontinuity Feedback}

\author[Zhang Z R et.~al.]{Hong Zhang\affil{1},Yue Zhao\affil{2},
       Xing Ji\affil{1}\comma\corrauth, and Kun Xu\affil{3,4,5}}
 \address{\affilnum{1}\ State Key Laboratory for Strength and Vibration of Mechanical Structures, School of Aerospace Engineering, Xi'an Jiaotong University, Xi'an 710000, P.R. China \\
 		   \affilnum{2}\ China Academy of Launch Vehicle Technology, Beijing, China\\
           \affilnum{3}\ Department of Mathematics, Hong Kong University of Science and Technology, Clear Water Bay, Kowloon 999077, Hong Kong SAR\\
           \affilnum{4}\ Department of Mechanical and Aerospace Engineering, Hong Kong University of Science and Technology, Clear Water Bay, Kowloon 999077, Hong Kong SAR\\
            \affilnum{5}\ Shenzhen Research Institute, Hong Kong University of Science and Technology, Shenzhen 518000, China}

 \emails{{\tt zhanghong2001@stu.xjtu.edu.cn} (H.~Zhang),\\ {\tt  kgfit@163.com} (Y.~Zhao), {\tt  jixing@xjtu.edu.cn} (X.~Ji), {\tt makxu@ust.hk} (K.~Xu)}

\begin{abstract}

This paper introduces an efficient class of adaptive stencil extension reconstruction methods based on a discontinuity feedback factor, addressing the challenges of weak robustness and high computational cost in high-order schemes, particularly those of 7th-order or above. Two key innovations are presented:
The accuracy order adaptively increases from the lowest level based on local stencil smoothness, contrasting with conventional methods like Weighted Essentially Non-Oscillatory (WENO) and Monotonic Upstream-Centered Scheme for Conservation Laws (MUSCL)limiters, which typically reduce order from the highest level.
The Discontinuity Feedback Factor (DF) serves a dual purpose: detecting sub-cell discontinuity strength and explicitly incorporating into the reconstruction process as a local smoothness measure. This approach eliminates the need for computationally expensive smoothness indicators often required in very high-order schemes, such as 9th-order schemes, and can be easily generalized to arbitrary high-order schemes.
Rigorous test cases, including a Mach 20000 jet, demonstrate the exceptional robustness of this approach.

\end{abstract}


\maketitle

\section{Introduction}
\label{sec1}

A primary focus in contemporary Computational Fluid Dynamics (CFD) is the development of high-order schemes for simulating turbulent flows in complex geometries. This research direction gained momentum following seminal work by Harten et al. \cite{harten1997uniformly}. Since then, prominent researchers have made significant contributions, resulting in several successful methodologies. Notable among these are:
Essential Non-Oscillatory (ENO) schemes \cite{harten1997uniformly,shu1988efficient}; Weighted Essential Non-Oscillatory (WENO) schemes \cite{jiang1996efficient,liu1994weighted}; Discontinuous Galerkin (DG) methods \cite{cockburn1989tvb,cockburn1998runge}.
These approaches have markedly advanced the field's capability to handle intricate flow problems with improved accuracy and efficiency.

High-order numerical methods face two main challenges: (1) The robustness of the algorithm will decrease rapidly as the order increases while calculate the discontinuity problems. Recently, there have been many strategies to improve the robustness of the algorithm. The mainstream approach is the new stencils selection strategy, such as the combination of the non-linear and linear weights \cite{gao2017enhanced}, which has the advantage of being able to improve the robustness while reduce the computational cost, such as combining two different WENO reconstructions and using a criterion to select the appropriate stencils \cite{liu2013robust}.
The TENO scheme provides new stencils  depending on the local flow characteristics, and can effectively handle situations with multiple discontinuities and recovers the robustness of the classical fifth-order WENO scheme \cite{jiang1996efficient,shu1988efficient}. Another effective way to improve the robustness is the limiter, such as the conventional prior limiters, like the van Leer limiter \cite{van1997towards}, the van Albada limiter \cite{van1982comparative} and post limiters like MOOD \cite{clain2011high}, etc. The parameters of the prior limiter usually have a significant effect on the results, and post limiters are usually computationally expensive. (2) High-order schemes are usually computationally expensive. For temporal advance, compared to classical Runge-Kutta method, the one-stage high-order method \cite{li2023one} can significantly reduces the computational cost. For spatial reconstruction, the WENO scheme requires an additional calculation of the smoothness indicators of the conservative variables to deal with the discontinuities in the stencil. The smoothness indicators will become complicated with significant increase of computational cost as the increase of the order and the number of the Gaussian points. Current approaches to reducing the computational cost focus on simplification of the smoothness indicators, such as the use of the combination of smoothness indicators of the lower-order stencils to linearly represent the smoothness indicators of the high-order stencil \cite{huang2018simple, kumar2018simple, kumar2019efficient}. Unfortunately, this method can not apply to arbitrary high-order schemes.
Another method is to use the weighted-least-squares method to give the smoothness indicators of each stencil \cite{huang2023multidimensional,liu2016wls}. 
However, the form of the global smoothness indicator $\tau_Z$ cannot be unified by the difference of the order, and the robustness of the method is weak.

In order to deal with the above challenges, the Adaptive Stencil Extension reconstruction methods based on Discontinuity Feedback factor (ASE-DF) are constructed.
As presented in Sec \ref{define of the discontinuity feedback factor}, DF is to define the discontinuity strength associated with the stencil-based reconstruction.  Without imposing the continuous flow distribution assumption inside each cell, the use of DF makes it more flexible in developing high-order FVM and its transition to be first-order scheme for preserving positivity property. The use of DF not only improves the robustness of the algorithm, 
but also replaces the smoothness indicators for each stencil, which significantly reduces the computational cost while achieving arbitrary high-order schemes. 
In the DF-based schemes presented in this paper, two flux functions are considered. The first one is the Lax-Friedrichs (L-F) with the SSP-RK \cite{Gottlieb2017EXPLICIT} temporal discretization, which has the positivity stability preserving, and the other  is the 2nd-order gas-kinetic scheme (GKS) \cite{xu1998gas,xu2001gks} with two-stages fourth-order (S2O4) temporal discretization.

This paper is organized as follows: Sec \ref{High-order finite volume scheme} provides a brief overview of the high-order finite volume scheme. Sec \ref{reconstruction method} introduces the ASP-DF reconstruction methods. In Sec \ref{numerical experiments}, some numerical results of the viscous and inviscid flow problems will be presented. The last section is the conclusion.

\section{High-order finite volume scheme}\label{High-order finite volume scheme}
\subsection{Finite volume framework}
The Finite Volume Method (FVM) involves dividing the computational domain into a finite amount of control volumes. Within each control volume, the physical quantities are integrated and averaged, resulting in the discrete equations
\begin{equation}
	\int_{\mathrm{\Omega}_{ij}}\mathbf{W}(\mathbf{x},t^{n+1}){\rm d}V=\int_{\mathrm{\Omega}_{ij}}\mathbf{W}(\mathbf{x},t^{n}){\rm d}V-\int^{t^{n+1}}_{t^n}\int_{\partial {\rm \Omega}_{ij}}\mathbf{F}\cdot \mathbf{n\rm d}S,
	\label{FVM discrete equations}
\end{equation}
where $\mathbf{W}=(\rho,\rho U, \rho V, \rho E)^T$ are the conservative variables in a control volume ${\rm \Omega}_{ij}$, and $\mathbf{F}$ are the fluxes across the $\partial {\rm \Omega}_{ij}$. $\partial {\rm \Omega}_{ij}$ corresponds to the cell interfaces, and in a 2-D rectangular mesh, the boundary $\partial {\rm \Omega}_{ij}$ can be expressed as
\begin{equation}
	\partial {\rm \Omega}_{ij}=\bigcup^4_{m=1}{\rm\Gamma}_{ij,m},
\end{equation}
where ${\rm \Gamma}_{ij,p}$ denotes the $p$th interface of the ${\rm \Omega}_{ij}$. Integrating over the cell ${\rm \Omega}_{ij}$, the semi-discrete form of the Eq (\ref{FVM discrete equations}) can be obtained as follows \cite{ji2020hweno}
\begin{equation}
	\frac{{\rm d}\mathbf{W}_{ij}}{{\rm d}t}=\mathcal{L}(\mathbf{W}_{ij})=-\frac1{|{\rm \Omega}_{ij}|}\sum^4_{m=1}\oint_{{\rm \Gamma}_{ij,m}}\mathbf{F}(\mathbf{W}_{ij})\cdot \mathbf{n}_m{\rm d}s,
\end{equation}
where $\mathbf{W}_{ij}$ is the cell average variables over the ${\rm\Omega}_{ij}$, $|{{\rm\Omega}_{ij}}|$ is the area of ${\rm\Omega}_{ij}$. $\mathcal{L}(\mathbf{W})$ is the temporal derivatives of the conservative variables, $\mathbf{F}=(F,G)^T$ are the flux function, and $\mathbf{n}_m$ corresponds to the outer normal direction of the interface ${\rm \Gamma}_{ij,m}$.

To achieve the high-order scheme, Gaussian points are considered. The line integral over ${\rm \Gamma}_{ij,m}$ is discretized according to Gaussian quadrature as follows
\begin{equation}
	\oint_{{\rm\Gamma}_{ij,m}}\mathbf{F}(\mathbf{W}_{ij})\cdot \mathbf{n}_m{\rm d}s=|l_m|\sum^p_{k=1}\eta_k\mathbf{F}(\mathbf{x}_{m,k},t)\cdot \mathbf{n}_m,
\end{equation}
where $\mathbf{x}_{m,k}$ for ${\rm \Gamma}_{ij,m}$ are the Gaussian points, $\eta_k$ is the weight of the $k$th Gaussian point, and $|l_m|$ corresponds to the length of the ${\Gamma}_{ij,m}$. For $r$th order spatial accuracy, $\frac{r-1}{2}$ Gaussian points are used.

To update the flow variables in global coordinates, first we need to obtain the conservative variables in the local coordinate
\begin{equation*}
	\widetilde{\mathbf{W}}=\mathbf{TW},
\end{equation*}
where the rotation matrix $\mathbf{T}$ for the 2-D case has the form
\begin{equation*}
	\mathbf{T}=
	\setlength{\arraycolsep}{8pt}
	\left(
	\begin{array}{cccc}
		1 & 0 & 0 & 0 \\
		0 & {\rm cos}\,\theta & {\rm sin}\,\theta & 0 \\
		0 & -{\rm sin}\,\theta & {\rm cos}\,\theta & 0 \\
		0 & 0 & 0 & 1 \\
	\end{array}
	\right).
\end{equation*}
Then we can obtain the fluxes in the local coordinate
\begin{equation}
	\widetilde{\mathbf{F}}(\mathbf{x}_{m,k},t)=\int \bm{\psi} \widetilde{u}f(\widetilde{\mathbf{x}}_{m,k},t,\widetilde{\mathbf{u}},\xi){\rm d}\widetilde{\mathbf{u}}{\rm d}\xi.
\end{equation}
According to \cite{2016Pan}, the global and local fluxes are related as
\begin{equation*}
	\mathbf{F}(\mathbf{x}_{m,k},t)\cdot \mathbf{n}=\mathbf{T}^{-1}\widetilde{\mathbf{F}}(\widetilde{\mathbf{x}}_{m,k},t),
\end{equation*}
here $\mathbf{F}$ and $\widetilde{\mathbf{F}}$ are the gas-kinetic fluxes, which will be presented in Sec \ref{gas flux function}.

\subsection{Gas-kinetic flux solver and two-stage fourth-order temporal discretization}\label{gas flux function}
The BGK equation \cite{Bhatnagar1954A} can be written as
\begin{equation}
	f_t+\mathbf{u}\cdot \nabla f=\frac{g-f}{\tau},
	\label{BGK equation}
\end{equation}
where $f$ is the gas distribution function, $g$ is the corresponding equilibrium state, and $\tau$ is the collision time. The collision term satisfied the compatibility condition
\begin{equation}
	\int \frac{g-f}{\tau}\bm{\psi}{\rm d}\Xi=0,
\end{equation}
where $\bm{\psi}=\left(1,u,v,\frac12(u^2+v^2+\xi^2)\right)^T$, ${\rm d}\Xi ={\rm d}u{\rm d}v{\rm d}\xi_1\cdots {\rm d}\xi_K$, $(u,v)$ are the two components of the macroscopic particle \color{black}velocities \color{black}, $\xi=(\xi_1,\cdots,\xi_K)$ correspond to the components \color{black}of \color{black} the internal particle \color{black}velocities \color{black} in $K$ dimensions, and $K$ is the number of the internal freedom. $K=(4-2\gamma)/(\gamma-1)$ for 2-D flows, and $\gamma$ is the specific heat ratio.

The general macroscopic gas dynamic equations can be obtained through the Chapman-Enskog expansion \cite{Chapman1970The} of the BGK equation, and the gas distribution function can be expressed as
\begin{equation}
	f=g-\tau D_{\mathbf{u}}g+\tau D_{\mathbf{u}}(\tau D_{\mathbf{u}})g-\tau D_{\mathbf{u}}[\tau D_{\mathbf{u}}(\tau D_{\mathbf{u}})g]+\cdots,
\end{equation}
where $D_{\mathbf{u}}=\partial/\partial t+\mathbf{u}\cdot\nabla$. When $f=g$, the equation takes the form of Euler equation, and the Navier-Stokes equation can be obtained by the truncated 1st-order distribution function
\begin{equation*}
	f=g-\tau(ug_x+vg_y+g_t).
\end{equation*}

The integral solution of the Eq (\ref{BGK equation}) is
\begin{equation}
	f(\mathbf{x}_{m,k},t,\mathbf{u},\xi)=\frac{1}{\tau}\int^t_0g(\mathbf{x}^\prime,t^\prime,\mathbf{u},\xi)e^{-(t-t^\prime)/\tau}dt^\prime+e^{-t/\tau}f_0(\mathbf{x}_{m,k},-\mathbf{u}t,\mathbf{u},\xi),
\end{equation}
where $\mathbf{x}_{m,k} = (0,0)$ is the quadrature point at the interface in the local coordinates, and $\mathbf{x}_{m,k}=\mathbf{x}^\prime+\mathbf{u}(t-t^\prime)$ is the trajectory of the particles, $f_0$ corresponds to the initial gas distribution function, and $g$ is the corresponding equilibrium state. Then we can construct a 2nd-order time accuracy gas distribution function \cite{xu2001gks} at the local Gaussian point $\mathbf{x}_{m,k}=(0,0)$
\begin{equation}
	\begin{aligned}
		&f(\mathbf{x}_{m,k},t,\mathbf{u},\xi)=\left(1-e^{-t/\tau_n}\right)g^c\\
		&+\left[(t+\tau)e^{-t/\tau_n}-\tau\right]a^c_\mathbf{x}\cdot\mathbf{u}g^c+\left(t-\tau+\tau e^{-t/\tau_n}\right)A^cg^c\\
		&+e^{-t/\tau_n}g^l\left[1-(\tau+t)a^l_\mathbf{x}\cdot\mathbf{u}-\tau A^l\right]\mathbb{H}(u)\\
		&+e^{-t/\tau_n}g^r\left[1-(\tau+t)a^r_\mathbf{x}\cdot\mathbf{u}-\tau A^r\right]\left[1-\mathbb{H}(u)\right],\\
	\end{aligned}
	\label{GKS flux}
\end{equation}
where $\mathbb{H}$ is the Heaviside function. The functions $g^{l,r,c}$ correspond to the initial gas distribution function on the left and right sides of the cell interface, and the equilibrium state located at the interface, respectively. The flow dynamics at the interface are contingent upon the ratio of the time step to the local particle collision time.

The function $g^k,k=l,r$ satisfies Maxwell's distribution
\begin{equation}
	g^k=\rho^k\left(\frac{\lambda^k}{\pi}\right)^{\frac{K+3}{2}}e^{-\lambda^k\left[(u-U^k)^2+(v-V^k)^2+\xi^2\right]},
\end{equation}
where $\lambda$ is a function of temperature, molecular mass and the Boltzmann constant. $g^k$ can be determined by
\begin{equation}
	\int\bm{\psi}g^l{\rm d}\Xi=\mathbf{W}^l,\quad \int\bm{\psi}g^r{\rm d}\Xi=\mathbf{W}^r,
	\label{relationship between g and W}
\end{equation}
where $\mathbf{W}^l,\mathbf{W}^r$ are the reconstructed variables, which will be presented in Sec \ref{reconstruction method}.

The coefficients $a_{\mathbf{x}},\ A$ denote the spatial and temporal derivatives, respectively, which have the form
\begin{equation*}
	a_{\mathbf{x}}\equiv(\partial g/\partial \mathbf{x})/g= g_{\mathbf{x}}/g,\quad A\equiv(\partial g/\partial t)/g= g_{t}/g,
\end{equation*}
which can be determined by the spatial derivatives of $\mathbf{W}$ and the compatibility condition as follows
\begin{equation*}
	\begin{aligned}
		&\left<a_x\right>=\frac{\partial \mathbf{W}}{\partial x}=\mathbf{W}_x,\quad \left<a_y\right>=\frac{\partial\mathbf{W}}{\partial y}=\mathbf{W}_{y},\\
		&\left<A+a_xu+a_yv\right>=0,
	\end{aligned}
\end{equation*}
where $a_x=(a_{x1},a_{x2}u,a_{x3}v,a_{x4}\frac12(u^2+v^2+\xi^2))^T$ and $a_y$ has the similar form. $\left<\cdots\right>$ are the moments of a gas distribution function defined by
\begin{equation}
	\rho\left<(\cdots)\right>=\int(\cdots)g{\rm d}\Xi.
\end{equation}
More details about the integration calculation can be found in Appendix \ref{GKS flux calculation}. Similarly, the equilibrium state $g^c$ and its derivatives $a^c_{\mathbf{x}}$ can be obtained by the reconstructed values $\mathbf{W}^c,\mathbf{W}^c_{\mathbf{x}}$, which will be presented in Sec \ref{center reconstruction}.

The physical collision time $\tau$ in the exponential function can be \color{black}modified \color{black} with a numerical collision time $\tau_n$ to capture unresolved discontinuities. For inviscid flow, the numerical collision time $\tau_n$ takes the following form \cite{xu2001gks}
\begin{equation}
	\tau_n=C_1\Delta t+C_2\left|\frac{p^l-p^r}{p^l+p^r}\right|\Delta t,
\end{equation}
where $C_1,C_2$ are constants, $p^l,p^r$ denote the pressure on the left and right side of the cell interface, respectively. For the \color{black}viscous \color{black} flow, $\tau_n$ is enlarged according to the normalized pressure difference \cite{xu2001gks}
\begin{equation}
	\tau_n=\tau+C_2\left|\frac{p^l-p^r}{p^l+p^r}\right|\Delta t=\frac{\mu}{p}+C_2\left|\frac{p^l-p^r}{p^l+p^r}\right|\Delta t,
\end{equation}
where $\mu$ is the dynamical viscosity coefficient. The reason for incorporating the pressure jump term to enlarge the particle collision time is to maintain the non-equilibrium dynamics within the shock layer via kinetic particle transport.

The two-stage fourth-order (S2O4) temporal discretization was developed in CFD applications \cite{LiangPan2016An,LiJiequan2015A}, which has the form
\begin{equation}\label{S2O4 temporal}
	\begin{aligned}
		&\mathbf{W}^*_{ij}=\mathbf{W}^n_{ij}+\frac12\Delta t\mathcal{L}(\mathbf{W}_{ij}^n)+\frac18\Delta t^2\frac{\partial}{\partial t}\mathcal{L}(\mathbf{W}_{ij}^n),\\
		&\mathbf{W}^{n+1}_{ij}=\mathbf{W}_{ij}^n+\Delta t\mathcal{L}(\mathbf{W}^n_{ij})+\frac{1}{6}\Delta t^2\left(\frac{\partial}{\partial t}\mathcal{L}(\mathbf{W}^n_{ij})+2\frac{\partial}{\partial t}\mathcal{L}(\mathbf{W}^*_{ij})\right),
	\end{aligned}
\end{equation}
where $\frac{\partial}{\partial t}\mathcal{L}(\mathbf{W})$ are the time derivatives of the flux transport integrated over the closed interfaces of the cell.

It is needed to obtain the first-order time derivatives of the flux at $t_n$ and $t^*=t_n+\Delta t/2$ \color{black}through the GKS flux function \color{black}. The total flux transport at the Gaussian point $x_{m,k}$ over the time interval $\delta$ is expressed as
\begin{equation}
	\mathbb{F}^n(\mathbf{x}_{m,k},\delta)=\int^{t_n+\delta}_{t_n}\mathbf{F}^n(\mathbf{x}_{m,k},t){\rm d}t=\int^{t_n+\delta}_{t_n}\int u\bm{\psi}f(\mathbf{x}_{m,k},t,\mathbf{u},\xi){\rm d}\Xi{\rm d}t.
	\label{total flux}
\end{equation}
For the 2nd-order GKS solver, let $t_n=0$, the flux can be approximated as a linear function
\begin{equation}
	\mathbf{F}^n(\mathbf{x}_{m,k},t)=\mathbf{F}^n(\mathbf{x}_{m,k})+\partial_t\mathbf{F}^n(\mathbf{x}_{m,k})t.
	\label{linear fulx}
\end{equation}
Substituting Eq (\ref{linear fulx}) into Eq.(\ref{total flux}), we can obtained
\begin{equation}
	\begin{aligned}
		&\mathbf{F}^n(\mathbf{x}_{m,k})\Delta t+\frac12\partial_t\mathbf{F}^n(\mathbf{x}_{m,k})\Delta t^2=\mathbb{F}^n(\mathbf{x}_{m,k},\Delta t),\\
		&\frac12\partial_t\mathbf{F}^n(\mathbf{x}_{m,k})\Delta t+\frac18\partial_t\mathbf{F}^n(\mathbf{x}_{m,k})\Delta t^2=\mathbb{F}^n(\mathbf{x}_{m,k},\Delta t/2).
	\end{aligned}
	\label{time advance}
\end{equation}
Eq (\ref{time advance}) is a linear system, and its solution is
\begin{equation*}
	\begin{aligned}
		&\mathbf{F}^n(\mathbf{x}_{m,k})=(4\mathbb{F}^n(\mathbf{x}_{m,k},\Delta t/2)-\mathbb{F}^n(\mathbf{x}_{m,k},\Delta t))/\Delta t,\\
		&\partial_t\mathbf{F}^n(\mathbf{x}_{m,k})=4(\mathbb{F}^n(\mathbf{x}_{m,k},\Delta t)-2\mathbb{F}^n(\mathbf{x}_{m,k},\Delta t/2))/\Delta t^2.
	\end{aligned}
\end{equation*}
Finally, the $\mathcal{L}(\mathbf{W}^n_{ij})$ and its time derivatives $\frac{\partial}{\partial t}\mathcal{L}(\mathbf{W}^n_{ij})$ can be obtained by
\begin{equation}
	\begin{aligned}
		&\mathcal{L}(\mathbf{W}^n_{ij})=-\frac1{{\Omega}_{ij}}\sum^4_{m=1}|l_m|\sum^p_{k=1}\eta_k\mathbf{F}^n(\mathbf{x}_{m,k})\cdot \mathbf{n}_m,\\
		&\frac{\partial}{\partial t}\mathcal{L}(\mathbf{W}^n_{ij})=-\frac1{{\Omega}_{ij}}\sum^4_{m=1}|l_m|\sum^p_{k=1}\eta_k\partial_t\mathbf{F}^n(\mathbf{x}_{m,k})\cdot \mathbf{n}_m.\\
	\end{aligned}
\end{equation}
In this way, all variables needed in Eq~(\ref{S2O4 temporal}) are determined.

\subsection{Lax-Friedrichs flux solver and RK temporal discretization}
To evaluate the robustness of the new reconstruction methods, the widely used Riemann solver, i.e., Lax-Freidrichs (L-F) flux is also considered. The 2-D Euler equations can be expressed as
\begin{equation}
	\left(
	\begin{array}{cccc}
		\rho \\
		\rho U \\
		\rho V \\
		\rho E \\
	\end{array}
	\right)_t+
	\left(
	\begin{array}{cccc}
		\rho U \\
		\rho U^2+p \\
		\rho UV \\
		U(\rho E+p) \\
	\end{array}
	\right)_x+
	\left(
	\begin{array}{cccc}
		\rho V \\
		\rho UV \\
		\rho V^2+p \\
		V(\rho E+p) \\
	\end{array}
	\right)_y=0.
\end{equation}
The L-F method \cite{Toro2013Riemann} is a numerical method for the solution of hyperbolic partial differential equations. The 1st-order L-F flux function can be expressed as
\begin{equation}
	\begin{aligned}
		\mathbf{F}(\mathbf{W}_{i+1/2})&=\frac12(\mathbf{F}(\mathbf{W}^l_{i+1/2})+\mathbf{F}(\mathbf{W}^r_{i+1/2}))\\
		&-\frac12{\rm max}\{|U^l_{i+1/2}|+c^l_{i+1/2},|U^r_{i+1/2}|+c^r_{i+1/2}\}(\mathbf{W}^r_{i+1/2}-\mathbf{W}^l_{i+1/2}),
	\end{aligned}
\end{equation}
\color{black}where $\mathbf{F}(\mathbf{W}^{l,r})=\left(\rho^{l,r}U^{l,r}, \rho^{l,r}U^{l,r}U^{l,r}+p^{l,r},\rho^{l,r}U^{l,r}V^{l,r},U^{l,r}(\rho^{l,r}E^{l,r}+p^{l,r})\right)^T$ correspond to the flux on the left and right side of the interface \color{black}, and $c^{l,r}=\sqrt{\frac{\gamma p^{l,r}}{\rho^{l,r}}}$ is the speed of sound.

The 3rd-order strong stability preserving Runge-Kutta (SSP-RK3) \cite{Gottlieb2017EXPLICIT} temporal discretization has been developed because of its strong \color{black}robust \color{black}, which has the form
\begin{equation}
	\begin{aligned}
		&\mathbf{W}^{(1)}_{ij}=\mathbf{W}^n_{ij}+\Delta t\mathcal{L}(\mathbf{W}_{ij}^n),\\
		&\mathbf{W}^{(2)}_{ij}=\frac34\mathbf{W}^n_{ij}+\frac14\mathbf{W}^{(1)}_{ij}+\frac14\Delta t\mathcal{L}(\mathbf{W}_{ij}^{(1)}),\\
		&\mathbf{W}^{n+1}_{ij}=\frac13\mathbf{W}^n_{ij}+\frac23\mathbf{W}^{(2)}_{ij}+\frac23\Delta t\mathcal{L}(\mathbf{W}_{ij}^{(2)}).\\
	\end{aligned}
	\label{SSP-RK3}
\end{equation}

\section{Adaptive stencil extension with discontinuous feedback reconstruction in zero-mean form}\label{reconstruction method}
In this section, an efficient class of high-order Adaptive Stencil Extension reconstructions with Discontinuous Feedback factor (ASE-DF) are presented. 
Compared to the classical WENO schemes, the ASE-DF schemes can improve computational efficiency, and significantly enhance the robustness of the algorithm.

\color{black}In this section, to reconstruct the left interface value $\mathbf{W}^l_{i+1/2,j}$ at the cell interface $x_{i+1/2,j}$, here we set $x_{i+1/2,j}=0$ for all the following reconstruct polynomials, and \color{black} for the domain $[-1,0]$, the zero-mean form polynomials are given by
\color{black}
\begin{align*}
	&Z_0(x)=1,\\
	&Z_1(x)=\frac{1}{\Delta x}x-\frac12,\quad Z_2(x)=\frac{1}{\Delta x^2}x^2+\frac13,\quad Z_3(x)=\frac{1}{\Delta x^3}x^3-\frac14,\quad Z_4(x)=\frac{1}{\Delta x^4}x^4+\frac15,\\
	&Z_5(x)=\frac{1}{\Delta x^5}x^5-\frac16,\quad Z_6(x)=\frac{1}{\Delta x^6}x^6+\frac17,\quad Z_7(x)=\frac{1}{\Delta x^7}x^7-\frac18,\quad Z_8(x)=\frac{1}{\Delta x^8}x^8+\frac19.
\end{align*}
\color{black}

\subsection{Discontinuity feedback factor}\label{define of the discontinuity feedback factor}
To deal with the possible discontinuities in the flow field, Ji et al. \cite{ji2021gradient, zhang2024robustness} proposed an indicator for feedback on the strength of the interface discontinuity based on the reconstructed values of the interface, which is called the discontinuity feedback factor (DF). The effect of DF is that when there exist discontinuities in the reconstruction stencil, the $n$th-order reconstruction polynomial will be automatically \color{black}reduced \color{black} to the 1st-order reconstruction, thus improving the robustness of the algorithm.

To capture the discontinuity more effectively, the improved DF is considered. First, we calculate the discontinuity strength $\sigma_{j+1/2}$ at the interface $\Gamma_{j+1/2}$
\begin{equation}
	\sigma_{j+1/2}={\mathbf{Avg}}\left\{\sum^n_{m=1}\sigma_{j+1/2,m}\right\},
\end{equation}
where $\sigma_{j+1/2}\ge0$, and $\sigma_{j+1/2,m}$ is the discontinuity strength of the $m$th Gauss point at the interface $\Gamma_{j+1/2}$, which has the form
\begin{equation}
	\begin{aligned}
		\sigma_{j+1/2,m}=\frac{\left|p^l-p^r\right|}{p^l}+\frac{\left|p^l-p^r\right|}{p^r}+\left({\rm Ma}^l_n-{\rm Ma}^r_n\right)^2+\left({\rm Ma}^l_t-{\rm Ma}^r_t\right)^2,
	\end{aligned}
\end{equation}
where $\rho^{k},k=l,r$ denote the left and right density of the Gauss point $x_{j+1/2.m}$, $p^{k},k=l,r$ denote the left and right pressure of the Gauss point $x_{j+1/2.m}$, ${\rm Ma}^l_n,\ {\rm Ma}^l_t$ denote the left Mach number defined by the normal and tangential velocity, respectively.

\paragraph{Remark 1} The $\sigma_{j+1/2}$ corresponds to the strength of the discontinuity. As long as the $\sigma_{j+1/2}=0$, the flow is smooth. However, the numerical solution itself is not strictly continuous, which means that some interfaces have a $\sigma_{j+1/2}$ of slightly larger than 0, \color{black}which corresponds to the weak discontinuities and can be dissipated during reconstruction, i.e. redundant reconstruction polynomial modifications. \color{black} To deal with the problem, we set a threshold $\sigma_{thres}$ for the $\sigma_{j+1/2}$, which is used to determine whether there is a discontinuity at the interface. The closer the threshold is to 0, \color{black}the more reconstruction polynomials will be modified by DF \color{black}. To improve the robustness and maintain the high resolution of the algorithm, a threshold $\sigma_{thres}$ is used in this paper, which means $\cdots+\sigma_{j-3/2}+\sigma_{j-1/2}+\sigma_{j+1/2}+\sigma_{j+3/2}+\cdots=0$ when $\cdots+\sigma_{j-3/2}+\sigma_{j-1/2}+\sigma_{j+1/2}+\sigma_{j+3/2}+\cdots<\sigma_{thres}$. The choice of $\sigma_{thres}$ will be analyzed in Sec \ref{discontinuity threshold}.
\\

Second, we need to calculate the DF factor $\alpha_{\mathbb{S}}$
\begin{equation}
	\begin{aligned}
		A&=\cdots+\sigma_{j-3/2}+\sigma_{j-1/2}+\sigma_{j+1/2}+\sigma_{j+3/2}+\cdots\\
		\alpha_{\mathbb{S}}&=\begin{cases}
			\quad 1.0& {\rm if\ } A<\sigma_{thres},\\
			\frac{\sigma_{thres}}{\cdots+\sigma_{j-3/2}+\sigma_{j-1/2}+\sigma_{j+1/2}+\sigma_{j+3/2}+\cdots}& {\rm otherwise },\\
		\end{cases}
	\end{aligned}
	\label{discontinuity strength}
\end{equation}
where $\alpha_{\mathbb{S}}\in(0,1]$ is the DF factor of the stencil $\mathbb{S}$. When $\alpha_{\mathbb{S}}= 1$, which means the stencil is smooth, when $\alpha_{\mathbb{S}}\rightarrow 0$, there are strong discontinuities in the stencil. $\{\cdots+\sigma_{j-3/2}+\sigma_{j-1/2}+\sigma_{j+1/2}+\sigma_{j+3/2}+\cdots\}$ is the sum discontinuity strength of all interfaces in a given direction in the stencil.

\begin{figure}[htbp]
	\centering
	\includegraphics[width=0.9\textwidth]{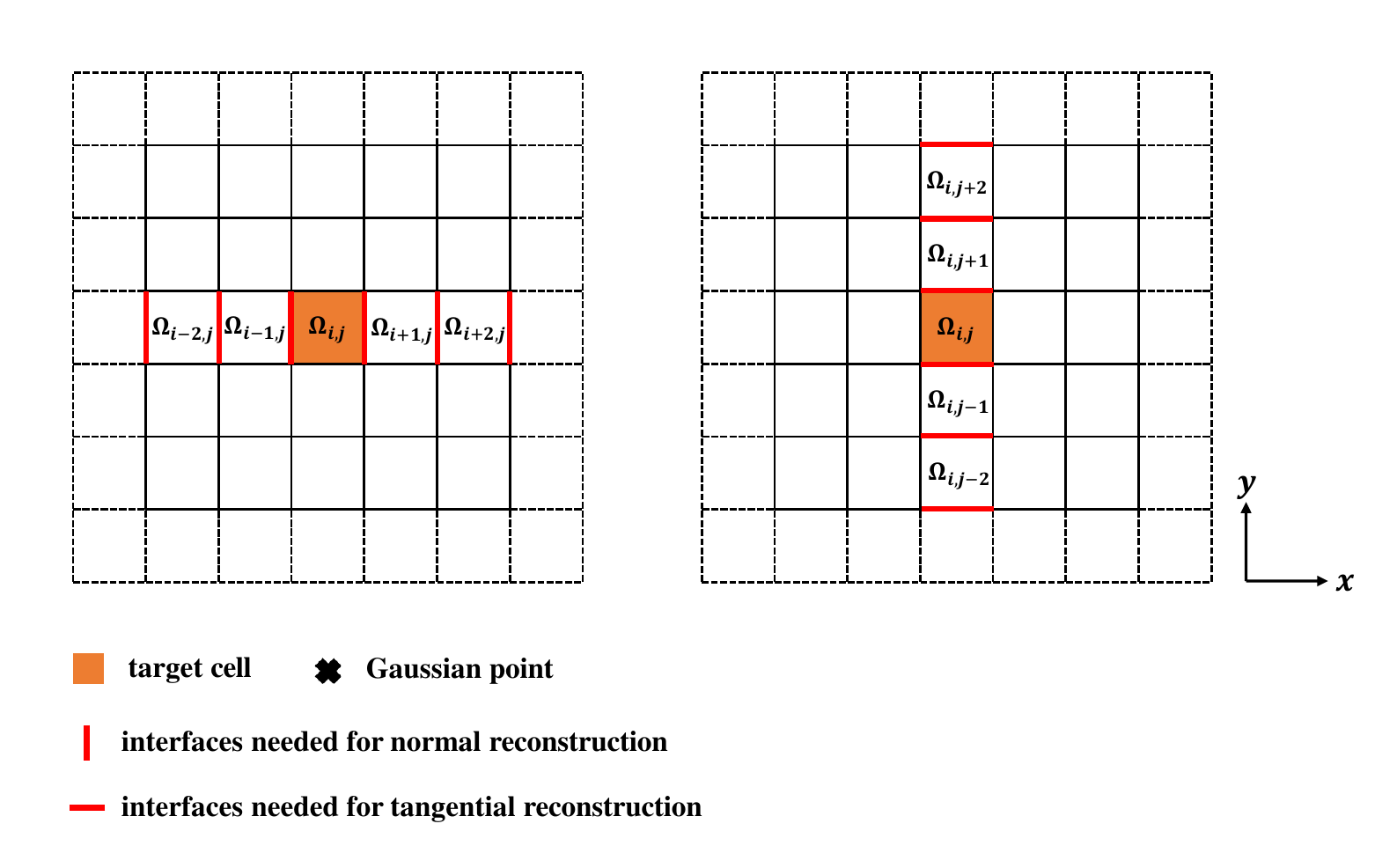}
	\caption{x-direction reconstruction at the left side of the interface ${\rm \Gamma}_{i+1/2,j}$. Discontinuity strength of the interfaces to be considered for normal and tangential reconstruction(from left to right).}
	\label{DFWENO schema}       
\end{figure}

\paragraph{Remark 2} For the normal and tangential reconstruction of the cell $\Omega_{i,j}$, it is necessary to calculate the DF of the corresponding direction of the stencil. As shown in Fig \ref{DFWENO schema}, taking the x-direction as an example, $\{\cdots+\sigma_{i-3/2,j}+\sigma_{i-1/2,j}+\sigma_{i+1/2,j}+\sigma_{i+3/2,j}+\cdots\}$ is considered in the normal reconstruction, and $\{\cdots+\sigma_{i,j-3/2}+\sigma_{i,j-1/2}+\sigma_{i,j+1/2}+\sigma_{i,j+3/2}+\cdots\}$ is considered in the tangential reconstruction.
\\

\subsection{ASE-DF(5,3)}

First we introduce the WENOZ-AO(5,3) method. The key idea of the 5th-order WENOZ-AO method is to reconstruct a reliable polynomial based on three 3rd-order sub-stencils $\{\mathbb{S}^{r3}_{-1},\mathbb{S}^{r3}_{0},\mathbb{S}^{r3}_{1}\}$ and a 5th-order stencil $\mathbb{S}^{r5}$ \cite{balsara2016efficient,borges2008improved}, which can obtain the 5th-order accuracy in the smooth region. Consider the variables $\{W_{-2},W_{-1}, W_0,W_1,W_2\}$, the stencil $\mathbb{S}^{r3}_{-1}$ gives
\begin{equation*}
	W_{x1}=2W_0-3W_{-1}+W_{-2},\quad W_{x2}=\frac12\left(W_{-2}-2W_{-1}+W_0\right).
\end{equation*}
The stencil $\mathbb{S}^{r3}_{0}$ gives
\begin{equation*}
	W_{x1}=W_1-W_0,\quad W_{x2}=\frac12\left(W_{-1}-2W_0+W_1\right).
\end{equation*}
The stencil $\mathbb{S}^{r3}_{1}$ gives
\begin{equation*}
	W_{x1}=W_1-W_0,\quad W_{x2}=\frac12\left(W_{0}-2W_{1}+W_2\right).
\end{equation*}
The stencil $\mathbb{S}^{r5}$ gives
\begin{equation*}
	\begin{aligned}
		W_{x1}&=\frac1{12}\left(-15W_0+W_{-1}+15W_1-W_2\right),\quad W_{x_2}= \frac18\left(-8W_0+6W_{-1}-W_{-2}+2W_1+W_2\right),\\
		W_{x3}&=\frac16\left(3W_0-W_{-1}-3W_1+W_2\right),\quad W_{x4}=\frac1{24}\left(6W_0-4W_{-1}+W_{-2}-4W_{1}+W_2\right).
	\end{aligned}
\end{equation*}

To deal with the discontinuity, the WENOZ type non-linear weights are used as
\begin{equation}
	\omega_k=d_k\left(1+\frac{\tau_Z^2}{(\beta_k+\varepsilon)^2}\right),
\end{equation}
where $\varepsilon=10^{-6}$, $d_k$ are the linear weights, and have the following values
\begin{equation*}
	d^5_0=d_{Hi},\ d^3_{-1}=(1-d_{Hi})(1-d_{Lo})/2,\ d^3_0=(1-d_{Hi})d_{Lo},\ d^3_1=d^3_{-1}.
\end{equation*}
Here $d_{Hi}=d_{Lo}=0.85$ are used. The global smoothness indicator $\tau_Z$ is defined as
\begin{equation*}
	\tau_Z=\frac13(|\beta^5_0-\beta^3_{-1}|+|\beta^5_0-\beta^3_1|+|\beta^5_0-\beta^3_1|).
\end{equation*}
The local smoothness indicators are defined as
\begin{equation}
	\beta_k=\sum^{q_k}_{q=1}\Delta x^{2q-1}\int^{x_{i+1/2,j}}_{x_{i-1/2,j}}\left(\frac{{\rm d}^q}{{\rm d}x^q}{\rm P}_k(x)\right)^2{\rm d}x,
	\label{LocalSmoothnessIndicator}
\end{equation}
where $q_k$ is the order of ${\rm P}_k(x)$. However, the calculation of the smoothness indicator is a heavy part of the whole reconstruction, for the WENOZ-AO reconstruction, not only $\beta^3_k,\ k=-1,0,1$ but also $\beta^5_0$ need to be calculated, and Huang et al. \cite{huang2018simple} found the computational cost of $\beta^5_0$ is comparable with the sum of one of all $\beta^3_k,\ k=-1,0,1$. To improve the computational efficiency of the algorithm, here we use a simple smoothness indicator $\overline{\beta}^5_0$ \cite{kumar2018simple} to replace the initial $\beta^5_0$ for ${\rm P}^{r5}(x)$
\begin{equation}
	\overline{\beta}^5_0=\frac16\left(\beta^3_{-1}+4\beta^3_{0}+\beta^3_1\right)+|\beta^3_{-1}-\beta^3_1|.
\end{equation}
\color{black}
\paragraph{Remark 3} To verify the consistency of the simplified $\overline{\beta^5_0}$ with the initial ${\beta^5_0}$, using the Taylor expansion of $W$ at $x_i$, we can obtain
\begin{equation*}
	\begin{cases}
		\beta^3_{-1}&=(W^{\prime})^2\Delta x^2+\left(\frac{13}{12}(W^{\prime\prime})^2-\frac23W^\prime W^{\prime\prime\prime}\right)\Delta x^4+O(\Delta x^5),\\
		\beta^3_{0}&=(W^{\prime})^2\Delta x^2+\left(\frac{13}{12}(W^{\prime\prime})^2+\frac13W^\prime W^{\prime\prime\prime}\right)\Delta x^4+O(\Delta x^5),\\
		\beta^3_{1}&=(W^{\prime})^2\Delta x^2+\left(\frac{13}{12}(W^{\prime\prime})^2-\frac23W^\prime W^{\prime\prime\prime}\right)\Delta x^4+O(\Delta x^5),\\
		\beta^5_0&=(W^{\prime})^2\Delta x^2+\frac{13}{12}(W^{\prime\prime})^2\Delta x^4+O(\Delta x^5),\\
	\end{cases}
\end{equation*}
whereas for smoothness indicator $\overline{\beta^5_0}$, we have
\begin{equation}
	\overline{\beta}^5_0=(W^{\prime})^2\Delta x^2+\frac{13}{12}(W^{\prime\prime})^2\Delta x^4+O(\Delta x^5),
\end{equation}
which is exactly same as $\beta^5_0$ upto truncation error level of order $O(\Delta x^5)$.
\\\\
\color{black}
Then, the global smoothness indicator $\tau_Z$ is modified as
\begin{equation*}
	\tau_Z=\frac13(|\overline{\beta}^5_0-\beta^3_{-1}|+|\overline{\beta}^5_0-\beta^3_1|+|\overline{\beta}^5_0-\beta^3_1|).
\end{equation*}

The ASE-DF method is to use DF factor to automatically select the appropriate stencils for reconstruction, and can be easily generalized to arbitrary high-order methods. The DF factor reflects the strength of the discontinuity in the stencil, and we consider the stencil to be smooth when $\alpha_{\mathbb{S}^{r5}}=1$, then we directly use the linear 5th-order polynomial reconstruction, which has the form

\begin{equation}
	\mathbb{P}_{DF}(x)=W_0+W_{x1}Z_1(x)+W_{x2}Z_2(x)+W_{x3}Z_3(x)+W_{x4}Z_4(x),
\end{equation}
when $\alpha_{\mathbb{S}^{r5}}<1$, it is necessary to fall back to the WENOZ-AO method, and to enhance the robustness of the algorithm, combined with the DF factor $\alpha_{\mathbb{S}^{r3}_i}, \alpha_{\mathbb{S}^{r5}}$, the $i$th reconstructed polynomial can be expressed as
\begin{equation*}
	\begin{aligned}
		\mathbb{P}^{r3}_i(x)&=W_0+\alpha_{\mathbb{S}^{r3}_i}\left[W_{x1}Z_1(x)+W_{x2}Z_2(x)\right],\\
		\mathbb{P}^{r5}(x)&=W_0+\alpha_{\mathbb{S}^{r5}}\left[W_{x1}Z_1(x)+W_{x2}Z_2(x)+W_{x3}Z_3(x)+Z_4(x)\right].\\
	\end{aligned}
\end{equation*}

Combined with the normalized weights $\overline{\omega}_k=\omega_k/\sum\omega_k$, the final form of the reconstructed polynomial is
\begin{equation}
	\mathbb{P}_{DF}(x)=\overline{\omega}^5_0\left(\frac1{d^5_0}\mathbb{P}^{r5}(x)-\sum^1_{k=-1}\frac{d^3_k}{d^5_0}\mathbb{P}^{r3}_k(x)\right)+\sum^1_{k=-1}\overline{\omega}^3_k\mathbb{P}^{r3}_k(x).
\end{equation}

\begin{figure}[h]
	\centering
	\includegraphics[width=0.8\textwidth]{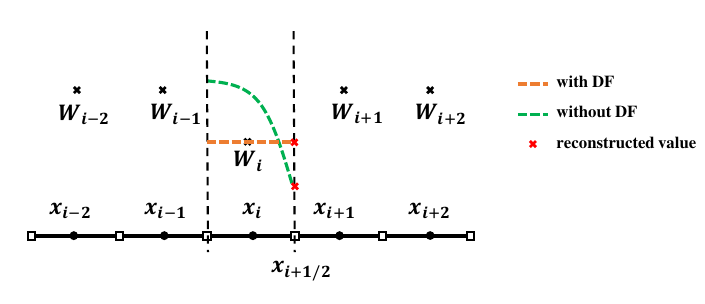}
	\caption{A possible distribution of variables $\{W_{i-2},W_{i-1},W_i,W_{i+1},W_{i+2}\}$. Each sub-stencil has a discontinuity, and the WENO reconstruction can only select the relatively smooth sub-stencil by weights. Take stencil $\{W_{i-2},W_{i-1},W_i\}$ as an example, the green line shows the reconstructed polynomial in the domain $[x_{i-1/2},x_{i+1/2}]$ by WENO method, and the orange line shows the DF factor can automatically converge the reconstructed polynomial to 1st-order when stencil exist a discontinuity, which is more robust compared to the WENO method.}
	\label{DF and WENO}       
\end{figure}

\paragraph{Remark 4}When $\alpha_{\mathbb{S}}\rightarrow1$, the reconstructed polynomial will revert to the initial smooth polynomial, when $\alpha_{\mathbb{S}}\rightarrow0$, which means there exists strong discontinuity in the stencil, and the reconstructed polynomial will be down to 1st-order, i.e. $\mathbb{P}(x)\approx W_0$. As shown in Fig \ref{DF and WENO}, $\alpha_{\mathbb{S}}$ can handle extreme cases that cannot be handled by the WENO scheme, thus improving the robustness of the algorithm.

\subsection{ASE-DF(7,5,3)}
The stencil $\mathbb{S}^{r7}$ gives
\begin{equation*}
	\begin{aligned}
		W_{x1}&=\frac1{180}(-245W_0 + 25W_{-1} - 2W_{-2} + 245W_1 - 25W_2 + 2f_3),\\
		W_{x2}&=\frac1{240}(-230W_0 + 210W_{-1} - 57W_{-2} + 7W_{-3} + 15W_1 + 63W_2 - 8W_3),\\
		W_{x3}&=\frac1{36}(28W_0 - 11W_{-1} + W_{-2}- 28W_1 + 11W_2 - W_3),\\
		W_{x4}&=\frac1{144}(46 W_0 - 39W_{-1} + 15W_{-2} - 2W_{-3} - 24W_1 + 3W_2 + W_3),\\
		W_{x5}&=\frac1{120}(-10 W_0 + 5 W_{-1} - W_{-2} + 10 W_1 - 5 W_2 + W_3),\\
		W_{x6}&=\frac1{720}(-20 W_0 + 15 W_{-1} - 6 W_{-2} + W_{-3} + 15 W_1 - 6 W_2 + W_3).
	\end{aligned}
\end{equation*}

When $\alpha_{\mathbb{S}^{r5}}=1$, which means that the stencil is smooth, we try to extend the two cells $\{W_{-3},W_3\}$ to get the stencil $\mathbb{S}^{r7}=\{W_{-3},W_{-2},W_{-1},W_{0},W_1,W_2,W_3\}$. Calculate the DF factor $\alpha_{\mathbb{S}^{r7}}$ of the stencil $\mathbb{S}^{r7}$ by Eq.(\ref{discontinuity strength}), when $\alpha_{\mathbb{S}^{r7}}=1$, which implies that $\mathbb{S}^{r7}$ is also smooth, then we directly use the linear 7th-order polynomial reconstruction, which has the form
\begin{equation}
	\begin{aligned}
		\mathbb{P}_{DF}(x)&=W_0+W_xZ_1(x)+Z_2(x)+Z_3(x)
		+W_{x4}Z_4(x)+W_{x5}Z_5(x)+W_{x6}Z_6(x),
	\end{aligned}
\end{equation}
when $\alpha_{\mathbb{S}^{r7}}<1$, the reconstruction will fall back to ASE-DF(5,3). Compared with the classical WENO method, the ASE-DF does not require additional calculation of the smoothness indicator, which significantly improves computational efficiency.

\subsection{ASE-DF(9,7,5,3)}
The stencil $\mathbb{S}^{r9}$ gives
\begin{align*}
	W_{x1}&=\frac{1}{5040}(-7175W_0 + 889W_{-1} - 119W_{-2} + 9W_{-3} + 7175W_1 - 889W_2 + 119 W_3- 9 W_4),\\
	W_{x2}&=\frac{1}{30240}(-27895W_0 + 28679W_{-1} - 9835W_{-2} + 2081W_{-3}- 205W_{-4} - 2065 W_1\\
	& + 11459W_2 - 2455W_3 + 236W_4),\\
	W_{x3}&=\frac1{1440}(1365W_0 - 587W_{-1} + 89W_{-2} - 7W_{-3} - 1365W_1 + 587W_2 - 89W_3 + 7W_4),\\
	W_{x4}&=\frac1{3456}(1174W_0 - 1160W_{-1} + 556W_{-2} - 128W_{-3} + 13 W_{-4} - 464W_1 - 68 W_2 + 88W_3 - 11 W_4),\\
	W_{x5}&=\frac1{480}(-75W_0 + 41W_{-1} - 11W_{-2} + W_{-3} + 75W_1 - 41W_2 + 11W_3 - W_4),\\
	W_{x6}&=\frac1{8640}(-380W_0 + 334W_{-1} - 170W_{-2}+ 46W_{-3} - 5W_{-4} + 250W_1 - 86W_2 + 10W_3 + W_4),\\
	W_{x7}&=\frac1{5040}(35W_0 - 21 W_{-1} + 7 W_{-2} - W_{-3} - 35 W_1 + 21 W_2 - 7W_3 + W_4),\\
	W_{x8}&=\frac1{40320}(70W_0 - 56W_{-1} + 28W_{-2} - 8W_{-3} + W_{-4} - 56W_1 + 28W_2 - 8 W_3 +W_4).
\end{align*}

Similarly to ASE-DF(7,5,3), when $\alpha_{\mathbb{S}^{r7}}=1$, we try to further expand the stencil and obtain the stencil $\mathbb{S}^{r9}$. Calculate the DF factor $\alpha_{\mathbb{S}^{r9}}$ by Eq (\ref{discontinuity strength}), when $\alpha_{\mathbb{S}^{r9}}=1$, then we directly use the linear 9th-order polynomial reconstruction, which has the form
\begin{equation}
	\begin{aligned}
		\mathbb{P}_{DF}(x)&=W_0+W_{x}Z_1(x)+W_{x2}Z_2(x)+W_{x3}Z_3(x)+W_{x4}Z_4(x)\\
		&+W_{x5}Z_5(x)+W_{x6}Z_6(x)+W_{x7}Z_7(x)+W_{x8}Z_8(x),
	\end{aligned}
\end{equation}
when $\alpha_{\mathbb{S}^{r9}}<1$, the reconstruction will fall back to ASE-DF(7,5,3). 

By expanding the stencil, we can easily obtain arbitrary higher-order reconstruction methods. The ASE-DF algorithm for selecting stencils is shown in Algorithm \ref{Alo1}.

\begin{algorithm}[h]\label{Alo1}
	\SetAlgoLined 
	\caption{DF-based adaptive stencil extension reconstruction methods}
	\KwIn{Interfaces reconstructed values $\mathbf{W}^l,\mathbf{W}^r$.}
	\KwOut{The needed stencil ${\rm \mathbb{S}}$ for the reconstruction.}
	\BlankLine
	calculate the correspond DF values $\alpha^{r5}\leftarrow(\mathbf{W}^l,\mathbf{W}^r)$ by Eq.(\ref{discontinuity strength}).\\
	\eIf{$\alpha^{r5}<1$}{
		select the non-linear WENOZ-AO(5,3) with DF reconstruction method;\\
	}{
	calculate the correspond DF values $\alpha^{r7}\leftarrow(\mathbf{W}^l,\mathbf{W}^r)$ by Eq.(\ref{discontinuity strength}).\\
		\eIf{$\alpha^{r7}<1$}{
			select the linear 5th-order polynomial reconstruction method;\\
		}{
		calculate the correspond DF values $\alpha^{r9}\leftarrow(\mathbf{W}^l,\mathbf{W}^r)$ by Eq.(\ref{discontinuity strength}).\\
			\eIf{$\alpha^{r9}<1$}{
				select the linear 7th-order polynomial reconstruction method;\\
			}{
			\color{black}
			calculate the correspond DF values $\alpha^{r11}\leftarrow(\mathbf{W}^l,\mathbf{W}^r)$ by Eq.(\ref{discontinuity strength}).\\
				\eIf{$\alpha^{r11}<1$}
				{
					select the linear 9th-order polynomial reconstruction method;\\
				}
				{
					$\cdots$
				}
				\color{black}
			}
		}
	}
	{return.}
\end{algorithm}

\subsection{Reconstruction of equilibrium state for gas-kinetic scheme}\label{center reconstruction}
For the non-equilibrium state, the reconstruction can be used to obtain the equilibrium state $g^c,g^c_x$ and $g^c_y$ directly by $g^k_{\mathbf{x}},k=l,r$, and a kinetic-based weighting method is used
\begin{equation}
	\begin{aligned}
		&\int \psi g^c{\rm d}\Xi=\mathbf{W}^c=\int_{u>0}\psi g^l{\rm d}\Xi+\int_{u<0}\psi g^r{\rm d}\Xi,\\
		&\int \psi g^c_{\mathbf{x}}{\rm d}\Xi=\mathbf{W}_{\mathbf{x}}^c=\int_{u>0}\psi g_{\mathbf{x}}^l{\rm d}\Xi+\int_{u<0}\psi g_{\mathbf{x}}^r{\rm d}\Xi.\\
	\end{aligned}
\end{equation}
In this way, all the variables needed in the algorithm have been determined.

\section{Numerical experiments}\label{numerical experiments}
In this section, ASE-DF reconstruction methods and two flux solvers are used to simulate several inviscid and viscous problems.

\subsection{Analysis of the discontinuity threshold $\sigma_{thres}$}
\label{discontinuity threshold}
In our proposed algorithm, the selection of the $\sigma_{thres}$ is crucial for the calculation of the strong shock wave problems. From Eq (\ref{discontinuity strength}), it can be seen that when a smaller value of $\sigma_{thres}$ is chosen, \color{black}the more reconstruction polynomials will be modified by DF \color{black}, which leads to \color{black}more order-reducing \color{black} of the reconstruction. \color{black}As a result \color{black}, the robustness of the whole algorithm is improved but the resolution will be significantly reduced. When a larger value of $\sigma_{thres}$ is chosen, which means that possible discontinuities in the stencil will not be recognized, and more higher-order smooth stencils will be used, the resolution of the algorithm will be improved, unfortunately, it leads to more pronounced oscillations and the robustness will decrease.

\begin{figure}[htbp]
	\centering
	\subfigure{
		\includegraphics[width=0.42\textwidth]{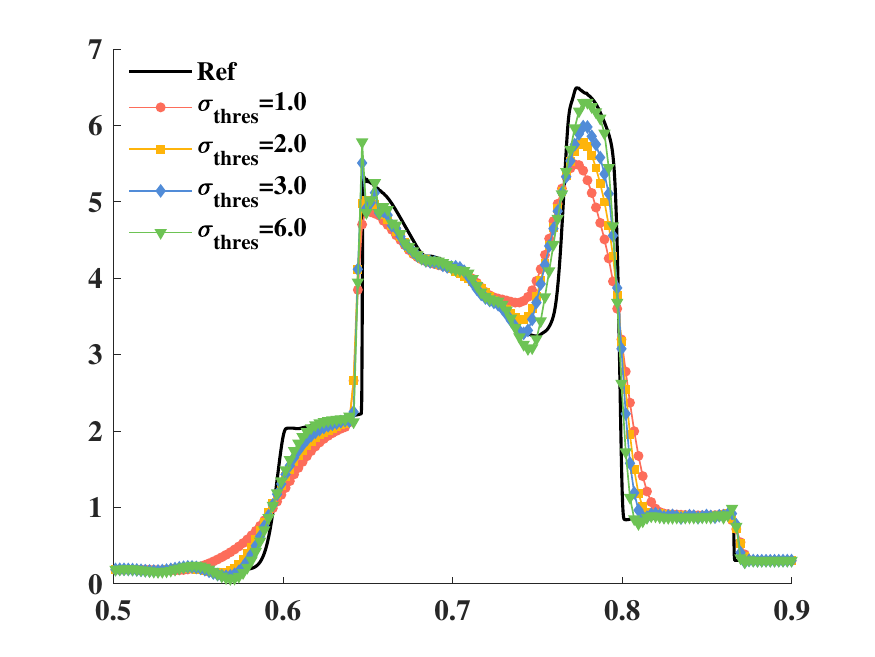}
	}
	\subfigure{
		\includegraphics[width=0.42\textwidth]{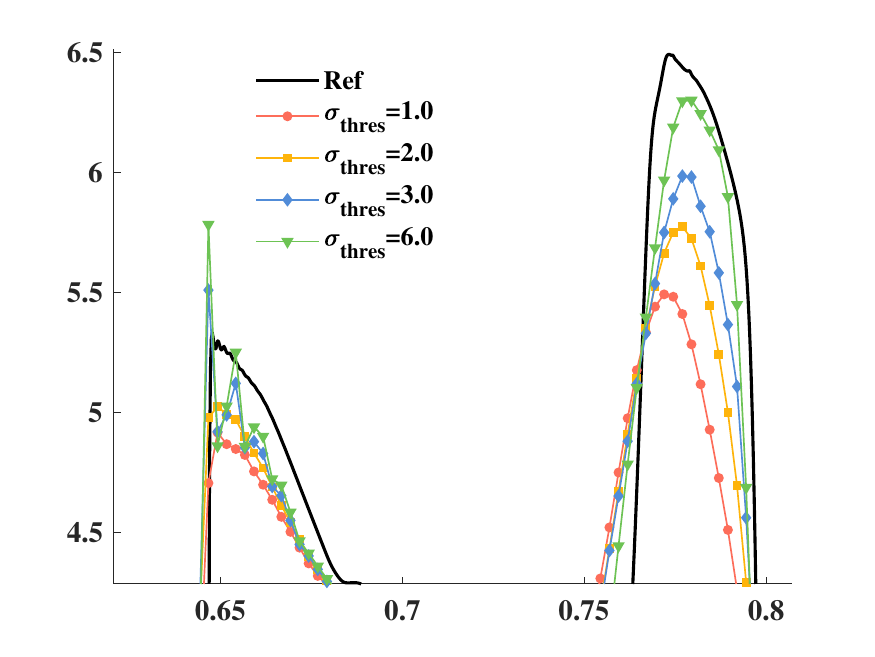}
	}
	\caption{Blast wave problem: the density distributions and local enlargement at $t=3.8$ with a cell size $\Delta x=1/400$. $c_1=0.05,c_2=5.0$. Different values of $\sigma_{thres}$ using the ASE-DF(5, 3) with the GKS solver. The reference solution is obtained by the 1-D 5th-order WENO-AO GKS with 4000 meshes.}
	\label{blast wave test}
\end{figure}

In order to select a suitable $\sigma_{thres}$, we present the numerical results of two strong shock wave problems.
\paragraph{Example 1}(Blast wave problem)\label{blast wave case test}
The initial conditions for the blast wave problem are given as follows
\begin{equation*}
	(\rho,u,p)=
	\begin{cases}
		(1,0,1000),&0\le x\le0.1,\\
		(1,0,0.01),&0.1\le x\le0.9,\\
		(1,0,100),&0.9\le x\le1.0.\\
	\end{cases}
\end{equation*}
400 uniform meshes are used and reflection boundary conditions are applied at both ends. The density distribution at $t=3.8$ is presented in Fig \ref{blast wave test}. It can be seen that as $\sigma_{thres}$ increases, the resolution of the algorithm increases with a significant increase in oscillations, while decreases significantly when $\sigma_{thres}$ is small.

\begin{figure}[htbp]
	\centering
	\subfigure[$\sigma_{thres}=1.0$]{
		\includegraphics[width=0.42\textwidth]{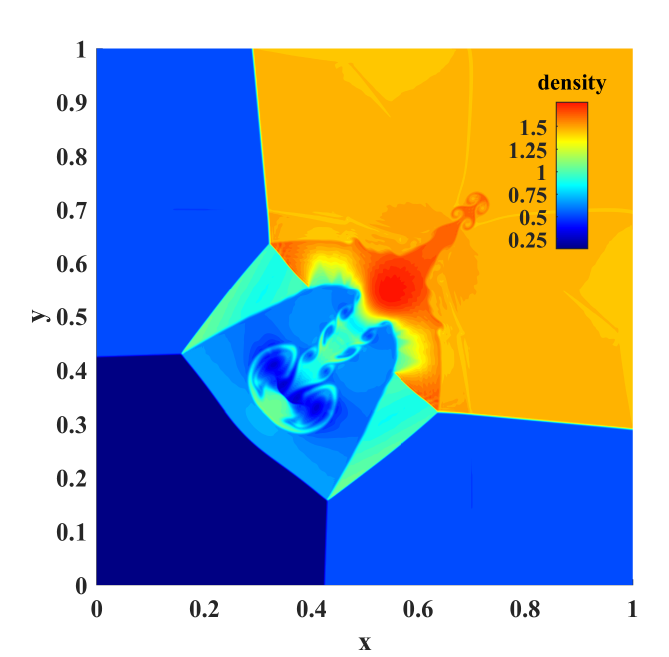}
	}
	\subfigure[$\sigma_{thres}=2.0$]{
		\includegraphics[width=0.42\textwidth]{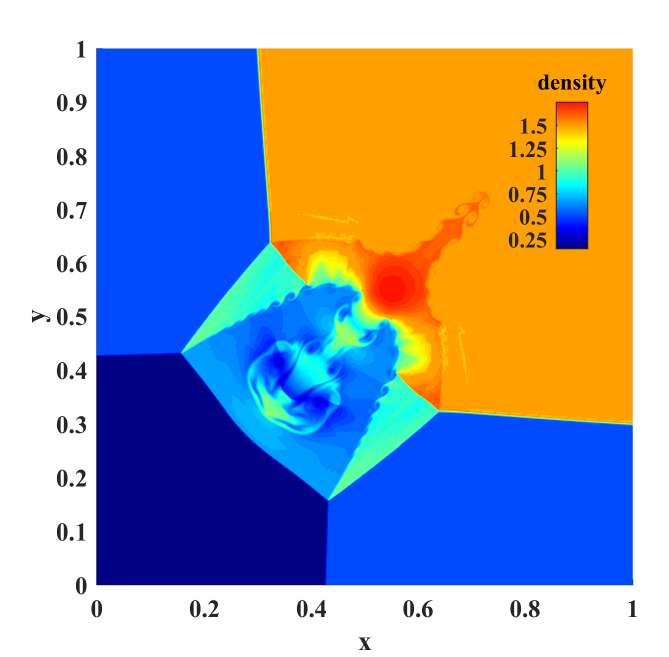}
	}
	\subfigure[$\sigma_{thres}=3.0$]{
		\includegraphics[width=0.42\textwidth]{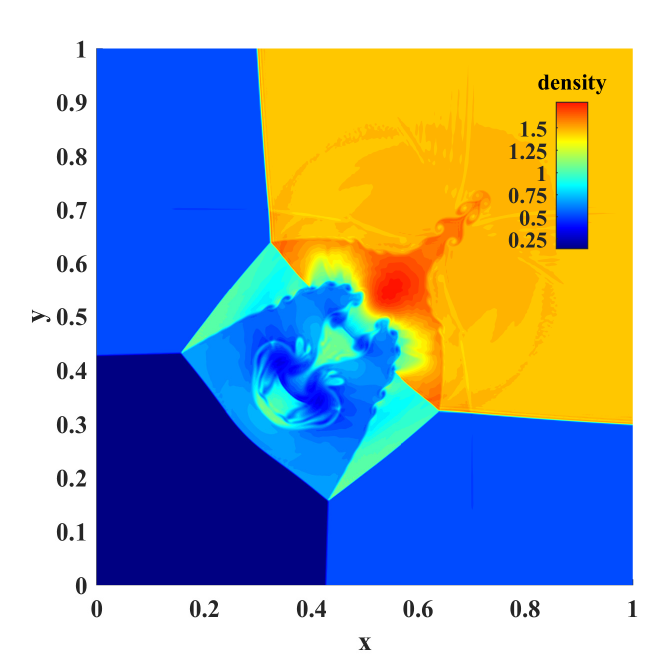}
	}
	\subfigure[$\sigma_{thres}=6.0$]{
		\includegraphics[width=0.42\textwidth]{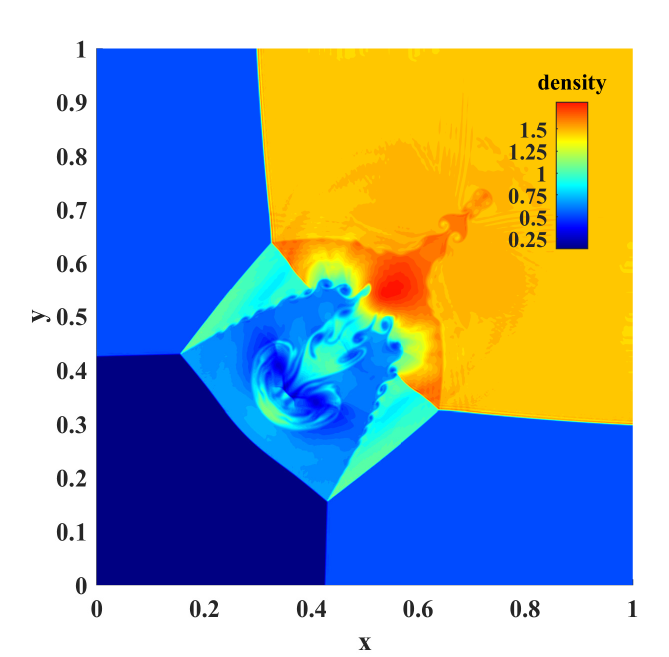}
	}
	\caption{Configuration 3: the density distribution at $t=0.6$ with $500\times500$ meshes. $c_1=0.05,c_2=1.0$. This figure is drawn with 30 density contours. (a-d) Different values of $\sigma_{thres}$ using the ASE-DF(5, 3) with the GKS solver.}
	\label{configuration 3 test}
\end{figure}

\paragraph{Example 2}(Configuration 3) Configuration 3 in \cite{Lax1998Solution} involves the shock-shock interaction and shock-vortex interaction. The initial condition in the domain $[0,1]\times[0,1]$ is given by
\begin{equation*}
	(\rho,u,v,p)=
	\begin{cases}
		(0.138,1.206,1.206,0.129),\quad &x<0.7,y<0.7,\\
		(0.5323,0,1.206,0.3),&x\ge 0.7,y<0.7,\\
		(1.5,0,0,1.5),&x\ge 0.7,y\ge 0.7,\\
		(0.5323,1.206,0,0.3),&x<0.7,y\ge 0.7.
	\end{cases}
\end{equation*}
The numerical results are shown in Fig \ref{configuration 3 test}, similar conclusion to blast wave case, further, no significant increase in resolution when $\sigma_{thres}$ is increasing.

Thus, in order to balance robustness and resolution of the algorithm, we select $\sigma_{thres}=2.0$ in the following numerical simulation.

\subsection{Accuracy validations}
To test the accuracy of the ASE-DF schemes, the smooth sin-wave propagation \cite{JiXing2021Compact} is considered. In these cases, both the physical viscosity and the collision time are set to zero. The initial condition of the equation is given by
\begin{equation}
	\rho(x)=1+0.2{\rm sin}(\pi x),\ U(x)=1.0,\ p(x)=1.0,\ x\in[0,2],
\end{equation}
and the exact solution has the form
\begin{equation*}
	\rho(x,t)=1.0+0.2{\rm sin}(\pi(x-t)),\ U(x,t)=1.0, \ p(x,t)=1.0.
\end{equation*}
Periodic boundary condition is used, and the numerical results are obtained after a periodic propagation at $t=2.0$. The accuracy of the schemes are measured in $L^1$-error norms in the domain $[0,2]$.

Extending to the 2-D case, the initial condition of the 2-D sin-wave propagation is
\begin{equation*}
	\begin{aligned}
		&\rho(x,y)=1.0+0.2{\rm sin}(\pi x){\rm sin}(\pi y),\\
		&U(x,y) = 1.0,\ V(x,y)=1.0,\ p(x,y)=1.0,
	\end{aligned}
\end{equation*}
with the exact solution
\begin{equation*}
	\begin{aligned}
		&\rho(x,y,t)=1.0+0.2{\rm sin}(\pi(x-t)){\rm sin}(\pi(y-t)),\\
		&U(x,y,t)=1.0,\ V(x,y,t)=1.0,\ p(x,y,t)=1.0.
	\end{aligned}
\end{equation*}

\begin{table}[htbp]
	\caption{Accuracy test for the 1-D sin-wave propagation in $L^1$-error norms.}
	\label{tab:1}       
	\begin{tabular*}{\textwidth}{@{\extracolsep{\fill}}ccccccc}
		\hline\noalign{\smallskip}
		N & ASP-DF(5,3) & order  & ASP-DF(7,5,3) & order & ASP-DF(9,7,5,3) & order\\
		\noalign{\smallskip}\hline\noalign{\smallskip}
		20 &2.585904e-05  &  &6.366042e-07 &  &1.424192e-08   &  \\
		40 &8.337553e-07  &4.95  &5.108345e-09 &6.96 &2.848478e-11  &8.97 \\
		80 &2.680325e-08  &4.95  &4.042055e-11 &6.98 &9.085371e-14 &8.29 \\
		160 &8.585536e-10   &4.96  &3.343673e-13 &6.92 &2.544534e-16 &8.48 \\
		\noalign{\smallskip}\hline
	\end{tabular*}
\end{table}
\begin{table}[htbp]
	\caption{Accuracy test for the 2-D sin-wave propagation in $L^1$-error norms.}
	\label{tab:2}       
	\begin{tabular*}{\textwidth}{@{\extracolsep{\fill}}ccccccc}
		\hline\noalign{\smallskip}
		N & ASP-DF(5,3) & order  & ASP-DF(7,5,3) & order & ASP-DF(9,7,5,3) & order\\
		\noalign{\smallskip}\hline\noalign{\smallskip}
		20$\times$20 &5.695852e-05  &  &2.475327e-06 &  &3.661497e-08  &  \\
		40$\times$40 &1.882335e-06  &4.92  &1.911540e-08  &7.02 &8.139979e-11  &8.81 \\
		80$\times$80 &6.106553e-08  &4.95  &1.490022e-10 &7.00 &4.792616e-13 &8.12 \\
		160$\times$160 &2.069689e-09  &4.88  &1.287922e-12 &6.85 &2.575173e-15 &7.54 \\
		\noalign{\smallskip}\hline
	\end{tabular*}
\end{table}

\color{black}To maintain the same accuracy in time and space, $dt=dx^{r/s}$ is used, where $r$ corresponds to the spatial order, and $s$ corresponds to the temporal order. \color{black} The numerical results are shown in Table \ref{tab:1}-\ref{tab:2}. The results show that the different schemes can achieve their theoretical numerical accuracy.

\subsection{Test cases with discontinuities}
\begin{figure}[htbp]
	\centering
	\subfigure{
		\includegraphics[width=0.42\textwidth]{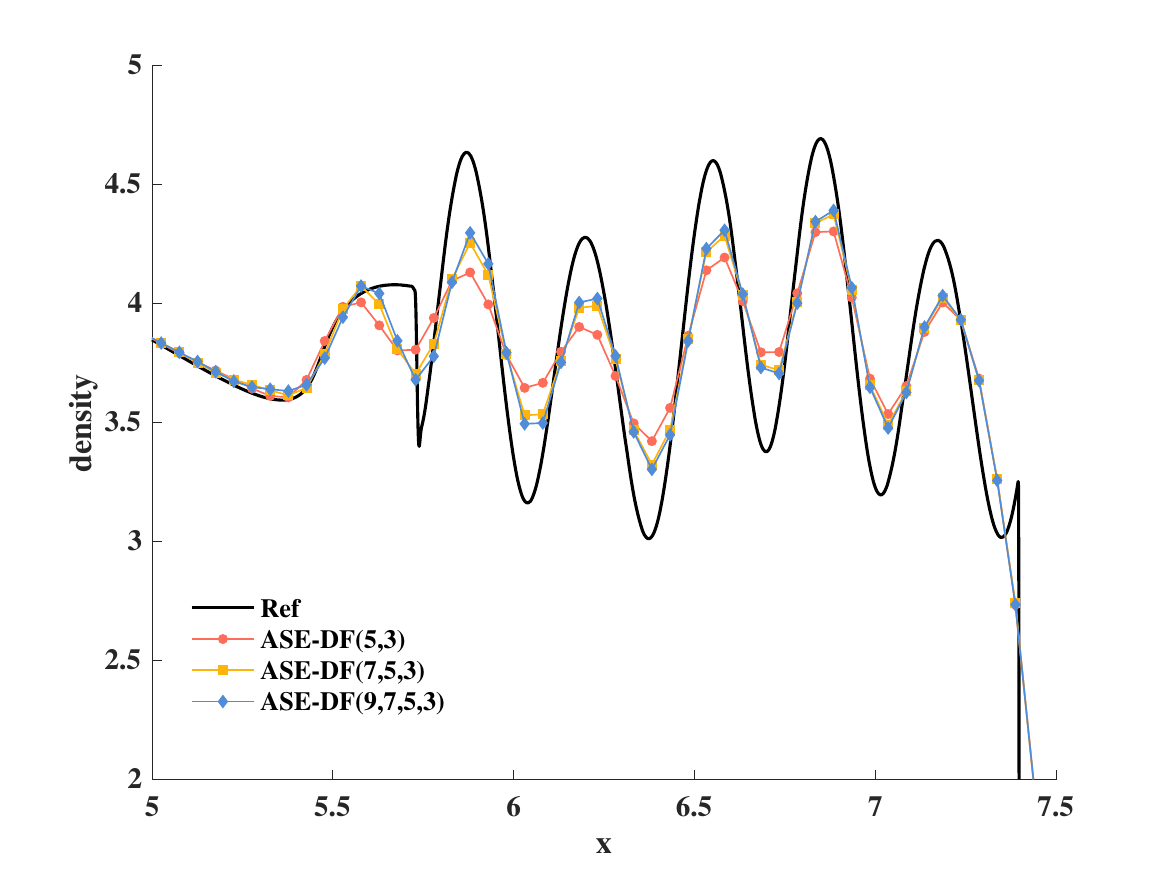}
	}
	\subfigure{
		\includegraphics[width=0.42\textwidth]{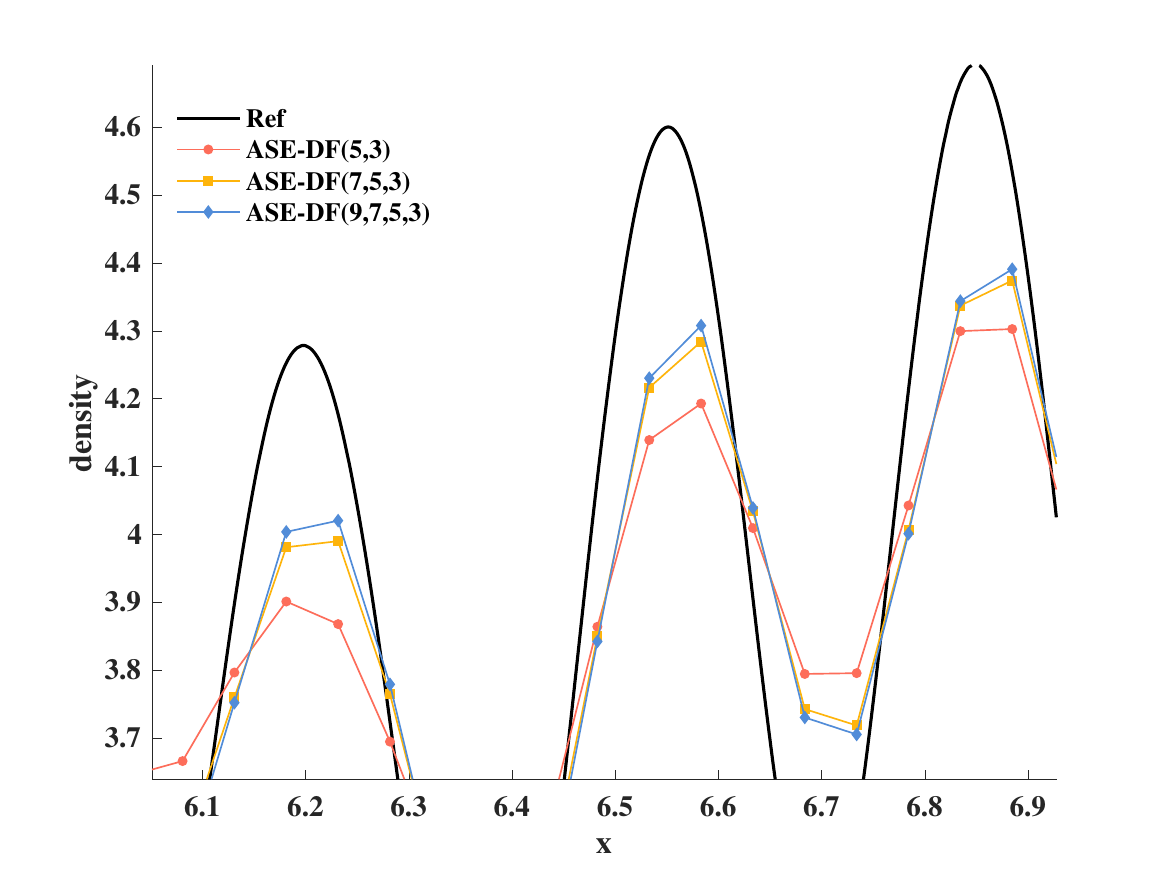}
	}
	\caption{Shu-Osher problem: the density distributions and local enlargement at $t=1.8$ with a cell size $\Delta x=1/20$. \color{black}All results are calculated by the GKS solver, $c_1=0.05,c_2=1.0$. \color{black} The reference solution is obtained by the 1-D 5th-order WENO-AO GKS with 10000 meshes.}
	\label{Shu-Osher}
\end{figure}

\paragraph{Example 1}(Shu-Osher problem) The Shu-Osher problem can examine the performance of capturing high frequency wave, and the initial condition is given by
\begin{equation*}
	(\rho,u,p)=
	\begin{cases}
		(3.857134, 2.629369, 10.33333),\quad &x\in[0.0,1.0],\\
		(1.0+0.2{\rm sin}(5x),0,1),\quad &x\in[1.0,10.0].
	\end{cases}
\end{equation*}
To clearly distinguish the ability of different schemes, 200 meshes are used. A CFL number $\theta_{CFL}=0.5$ is used for the all following cases. Fig \ref{Shu-Osher} shows the density distribution using the different reconstruction methods, for the same mesh size, both ASE-DF(7,5,3) and ASE-DF(9,7,5,3) perform better in in capturing high frequency waves compared to ASE-DF(5,3), which is consistent with theoretical expectations, and limited by the number of meshes, the resolution of ASE-DF(9,7,5,3) is not significantly higher compared to ASE-DF(7,5,3).

\begin{figure}[htbp]
	\centering
	\subfigure{
		\includegraphics[width=0.44\textwidth]{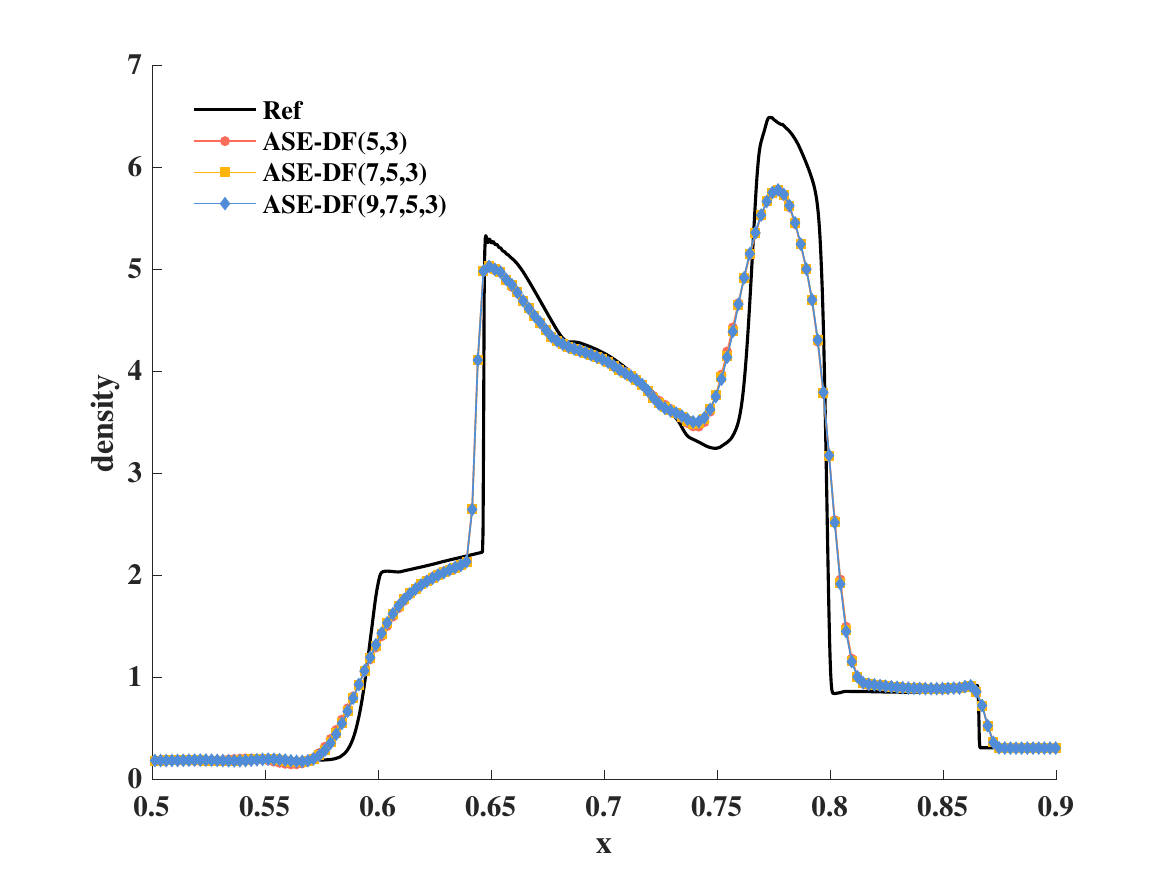}
	}
	\subfigure{
		\includegraphics[width=0.44\textwidth]{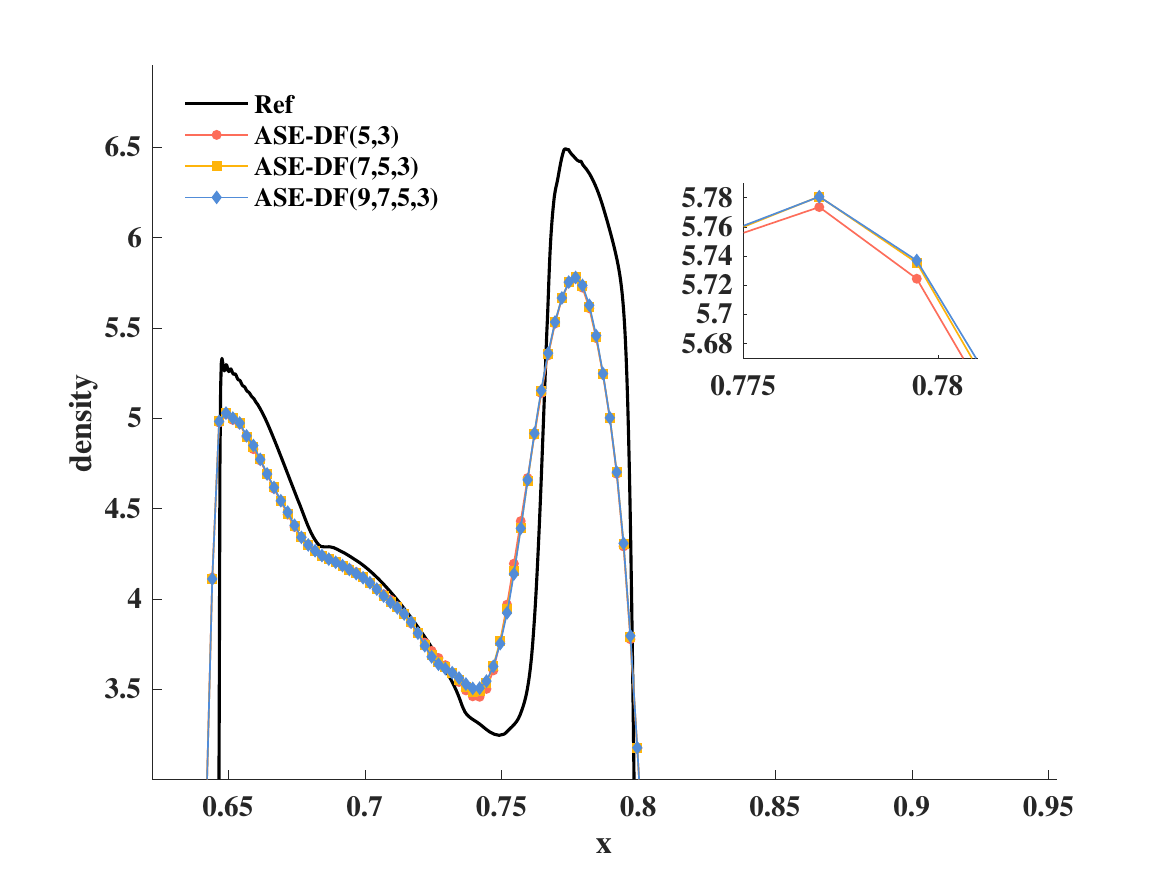}
	}
	\caption{Blast wave problem: the density distributions and local enlargement at $t=3.8$ with a cell size $\Delta x=1/400$. \color{black}All results are calculated by the GKS solver, $c_1=0.05,c_2=5.0$. \color{black} The reference solution is obtained by the 1-D 5th-order WENO-AO GKS with 4000 meshes.}
	\label{blast wave}
\end{figure}

\begin{figure}[h]
	\centering
	\subfigure[ASE-DF(5,3), GKS solver]{
		\includegraphics[width=0.31\textwidth]{figure/fig7.pdf}
	}
	\subfigure[ASE-DF(7,5,3), GKS solver]{
		\includegraphics[width=0.31\textwidth]{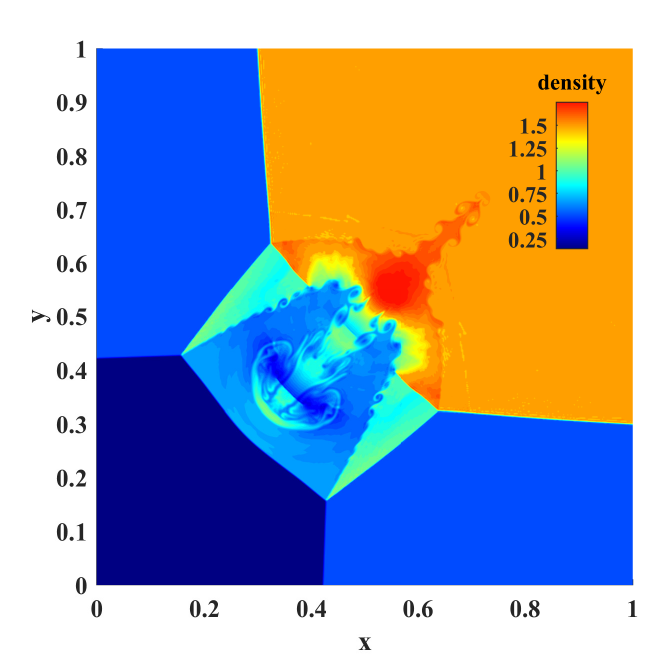}
	}
	\subfigure[ASE-DF(9,7,5,3), GKS solver]{
		\includegraphics[width=0.31\textwidth]{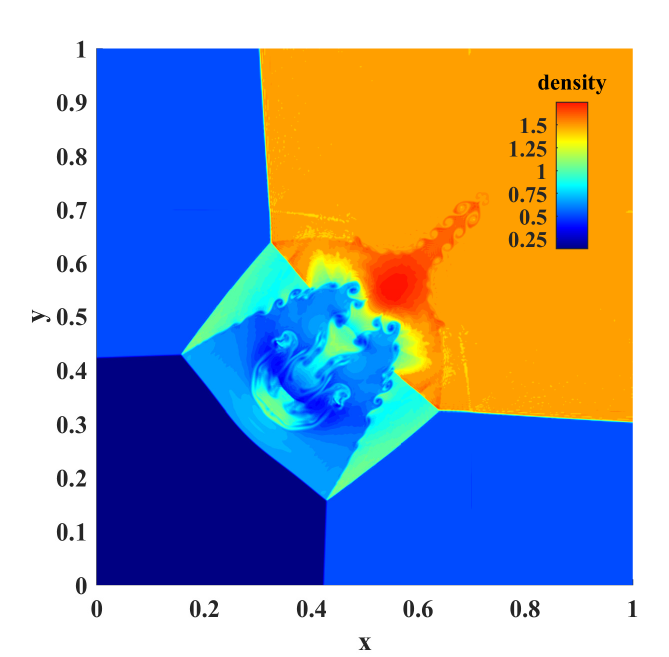}
	}
	\subfigure[ASE-DF(5,3), L-F solver]{
		\includegraphics[width=0.31\textwidth]{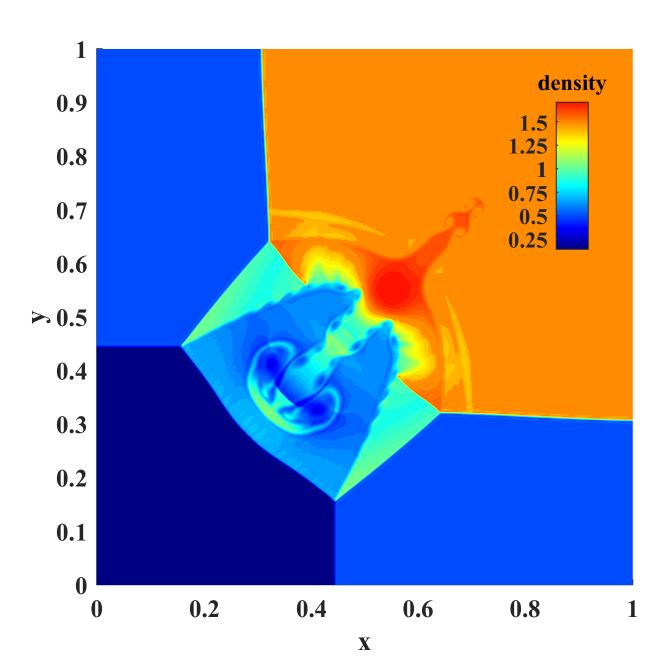}
	}
	\subfigure[ASE-DF(7,5,3), L-F solver]{
		\includegraphics[width=0.31\textwidth]{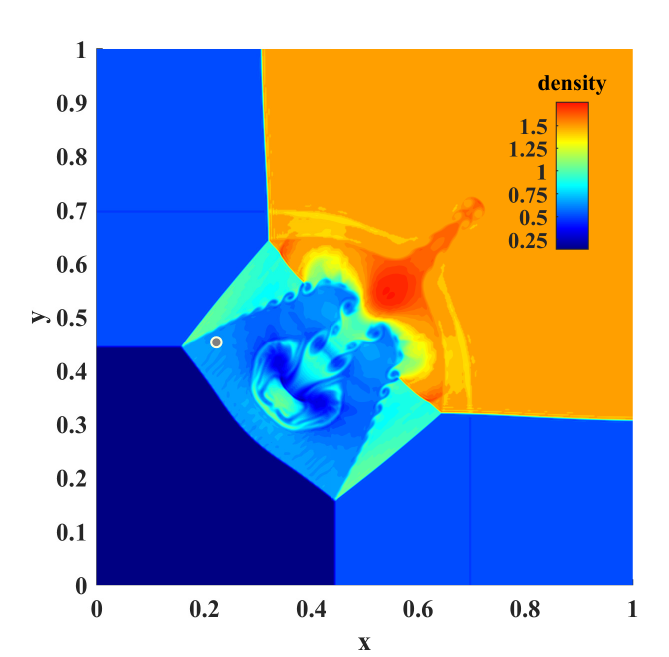}
	}
	\subfigure[ASE-DF(9,7,5,3), L-F solver]{
		\includegraphics[width=0.31\textwidth]{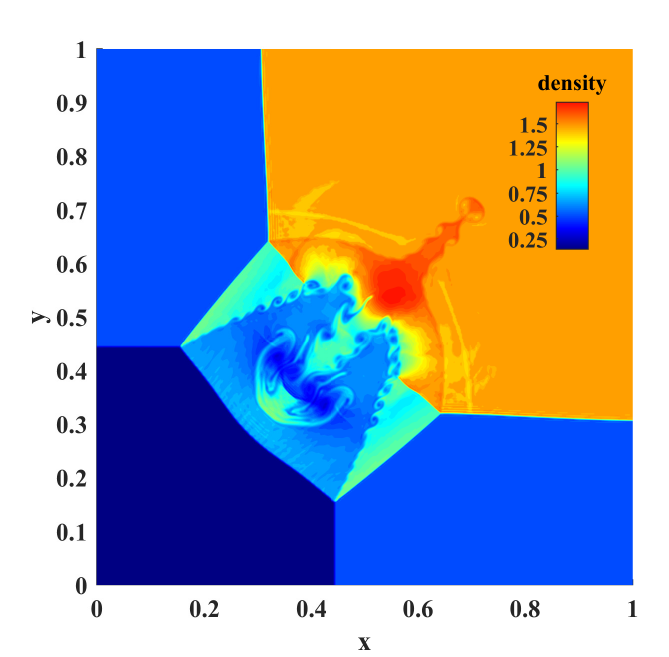}
	}
	\caption{Configuration 3: the density distribution at $t=0.6$ with $500\times500$ meshes. This figure is drawn with 30 density contours. (a-c) \color{black}The results are calculated by the GKS solver. $c_1=0.05,c_2=1.0$. (d-f) The results are calculated by the L-F solver. \color{black}}
	\label{configuration 3}
\end{figure}

\paragraph{Example 2}(Blast wave problem)\label{blast wave case}
To test the performance of the different reconstruction methods, the blast problem is simulated again. 400 uniform meshes are used and reflection boundary conditions are applied at both ends. The density distribution at $t=3.8$ is presented in Fig \ref{blast wave}. It can be seen that resolution of the algorithm is slightly increasing as the order increases, this is due to the fact that the presence of strong shock waves in this problem, and due to the DF, leads to higher-order reconstruction being treated more significantly by downgrading, i.e., ASE-DF(9,7,5,3) also perform almost identically to ASE-DF(5,3).

\begin{figure}[h]
	\centering
	\subfigure[ASE-DF(5,3), GKS solver]{
		\includegraphics[width=0.3\textwidth]{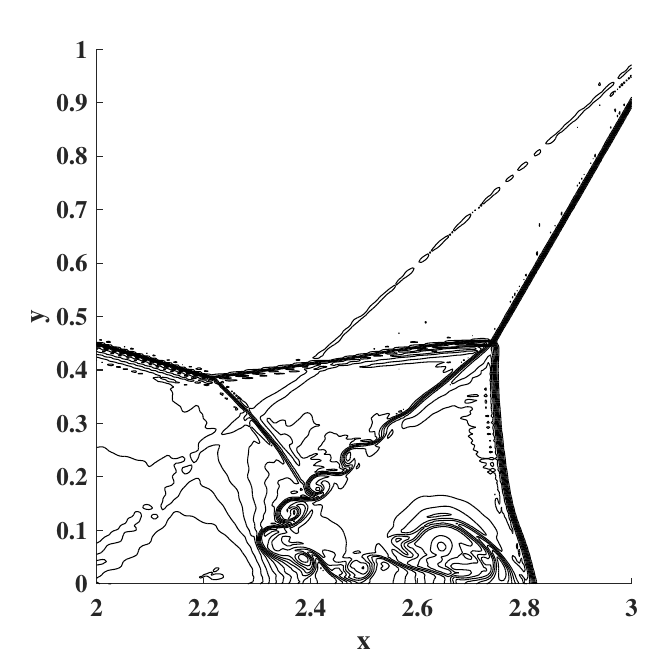}
	}
	\subfigure[ASE-DF(7,5,3), GKS solver]{
		\includegraphics[width=0.3\textwidth]{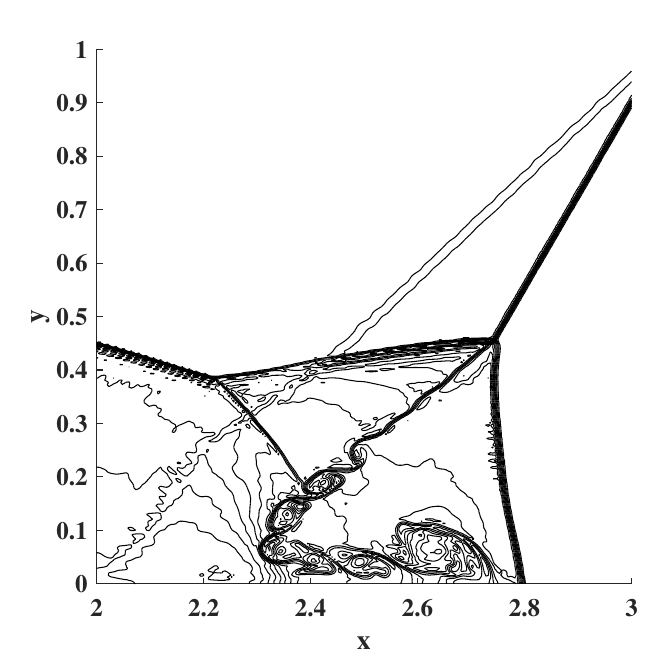}
	}
	\subfigure[ASE-DF(9,7,5,3), GKS solver]{
		\includegraphics[width=0.3\textwidth]{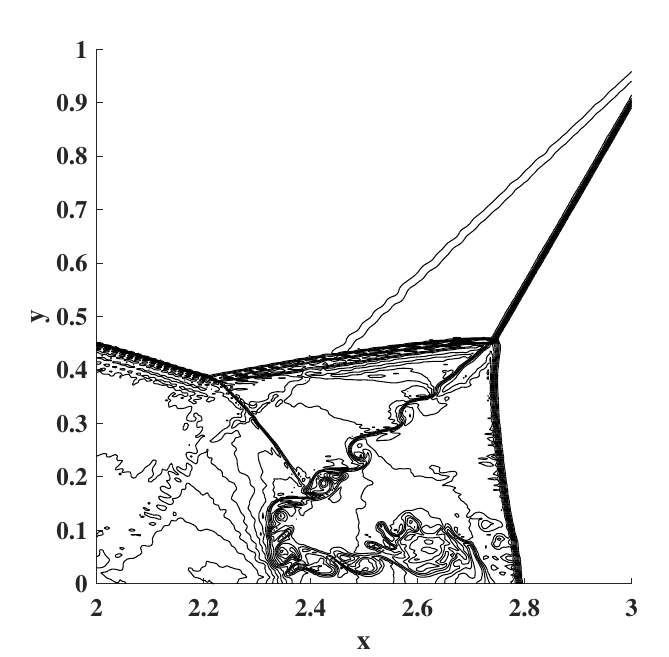}
	}
	\subfigure[ASE-DF(5,3), L-F solver]{
		\includegraphics[width=0.3\textwidth]{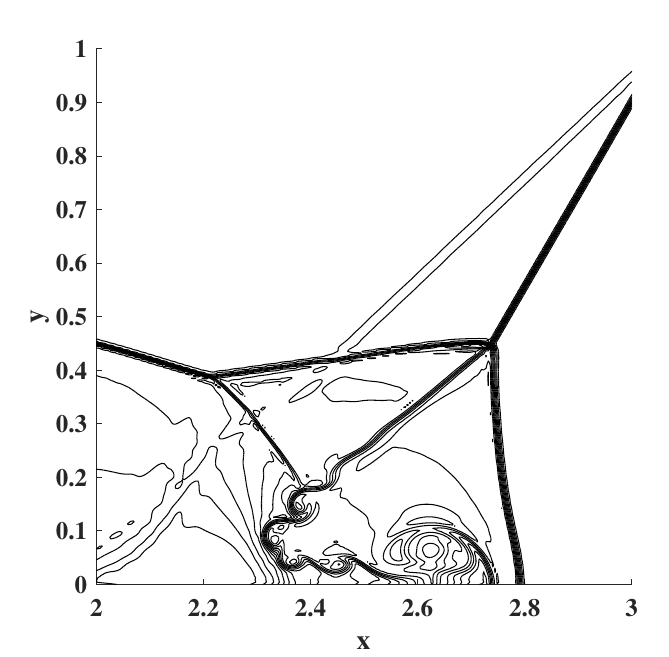}
	}
	\subfigure[ASE-DF(7,5,3), L-F solver]{
		\includegraphics[width=0.3\textwidth]{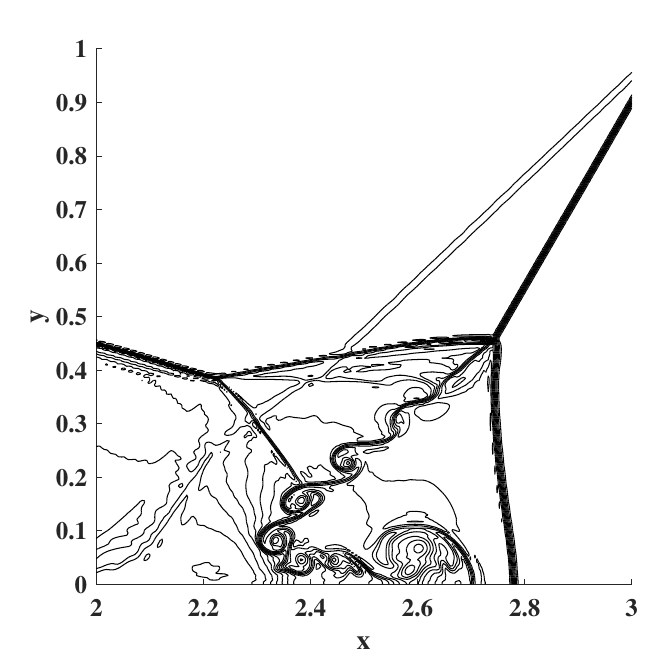}
	}
	\subfigure[ASE-DF(9,7,5,3), L-F solver]{
		\includegraphics[width=0.3\textwidth]{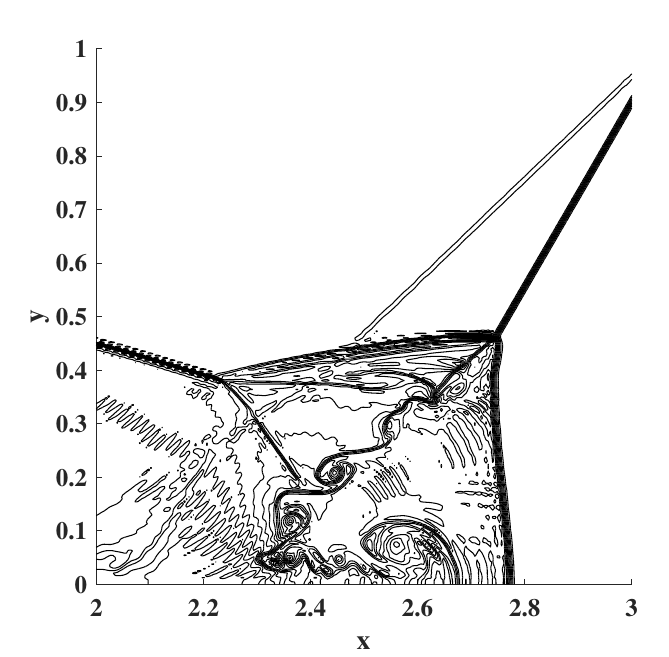}
	}
	\caption{double Mach reflection problem: the density distribution at $t=0.2$ with $960\times240$ meshes. This figure is drawn with 35 density contours. (a-c) \color{black}The results are calculated by the GKS solver. $c_1=0.05,c_2=1.0$. (d-f) The results are calculated by the L-F solver. \color{black}}
	\label{double Mach reflection}
\end{figure}

\paragraph{Example 3}(Configuration 3) The problem has been described in Sec \ref{discontinuity threshold}.
The numerical results are shown in Fig \ref{configuration 3}, all schemes show low-dissipation, and the main differences among all schemes are the strength of the shear layers, both ASE-DF(9,7,5,3) and ASE-DF(7,5,3) are capable of resolving significantly more vortical structures compared to ASE-DF(5,3).
\paragraph{Example 4}(double Mach reflection problem) In the double Mach reflection problem \cite{Woodward1984}, a rightwards-moving shock wave at Mach 10 with an incident angle of $60^{\circ}$ with respect to the x-axis interacting with a reflection wall boundary, which can test the shock wave capture capability, numerical stability, dissipation and resolution of the scheme. The initial condition is given by
\begin{equation*}
	(\rho,u,v,p)=
	\begin{cases}
		(1.4,0,0,1),\quad &y<1.732(x-0.1667),\\
		(8,7.145,-4.125,116.8333),&{\rm ohterwise.}\\
	\end{cases}
\end{equation*}
The slip boundary condition is used at the wall starting from $x=1/6$, the bottom boundary condition is the post-shock condition. The computational domain is $[0,4]\times[0,1]$, and the numerical results in $[2,3]\times[0,1]$ for all schemes at $t=0.2$ are shown in Fig \ref{double Mach reflection}, which can be observed that all schemes show low-dissipation and can capture the shock wave well.

\begin{figure}[htbp]
	\centering
	\subfigure[ASE-DF(5,3), $\Delta x=\Delta y=1/500$]{
		\includegraphics[width=0.42\textwidth]{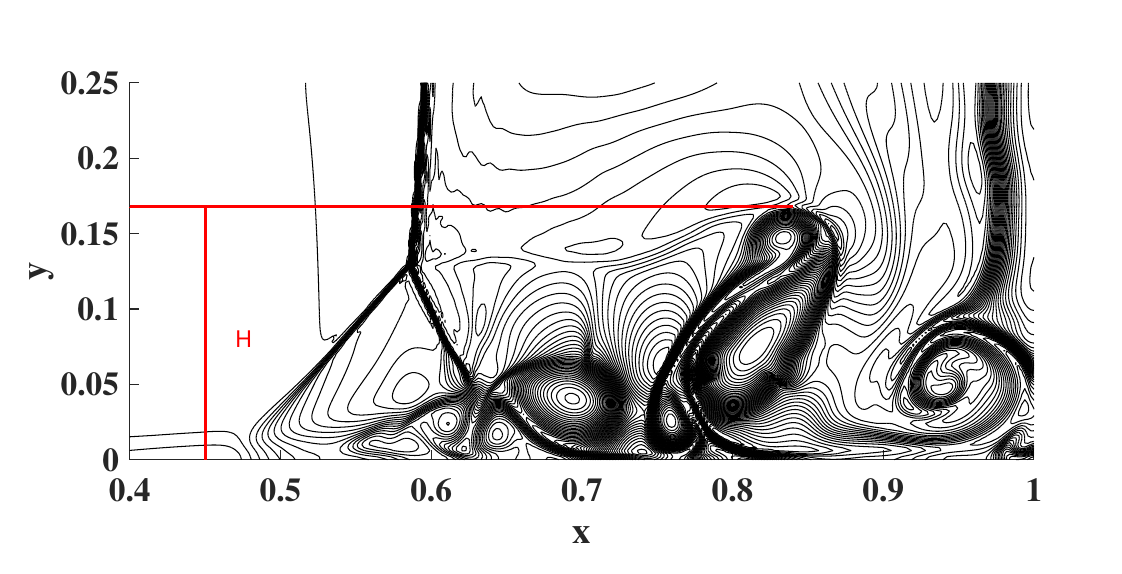}
	}
	\subfigure[ASE-DF(5,3), $\Delta x=\Delta y=1/1000$]{
		\includegraphics[width=0.42\textwidth]{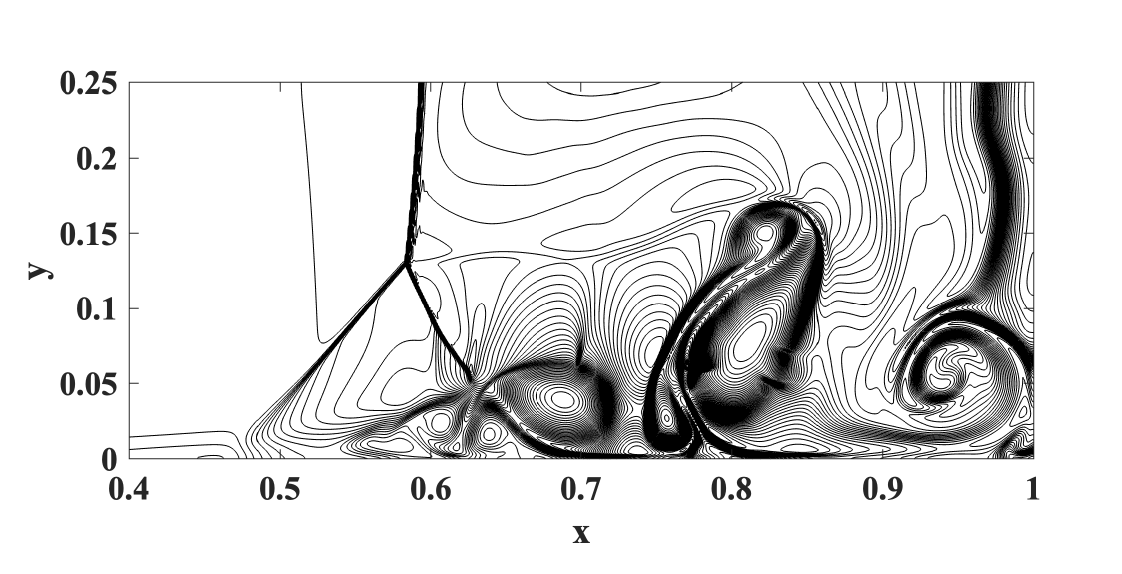}
	}
	\subfigure[ASE-DF(7,5,3), $\Delta x=\Delta y=1/500$]{
		\includegraphics[width=0.42\textwidth]{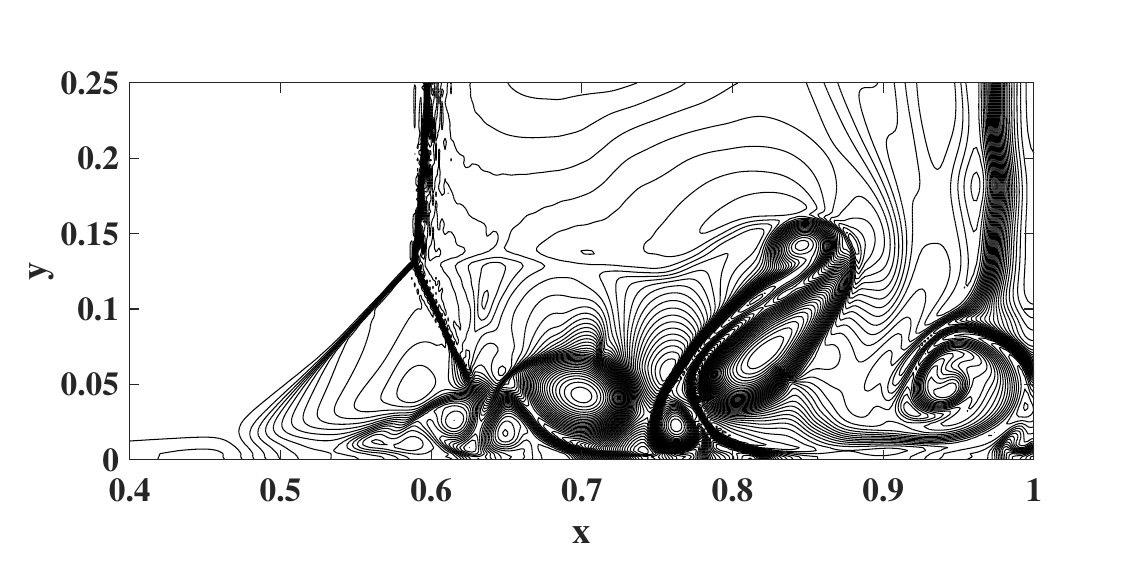}
	}
	\subfigure[ASE-DF(7,5,3), $\Delta x=\Delta y=1/1000$]{
		\includegraphics[width=0.42\textwidth]{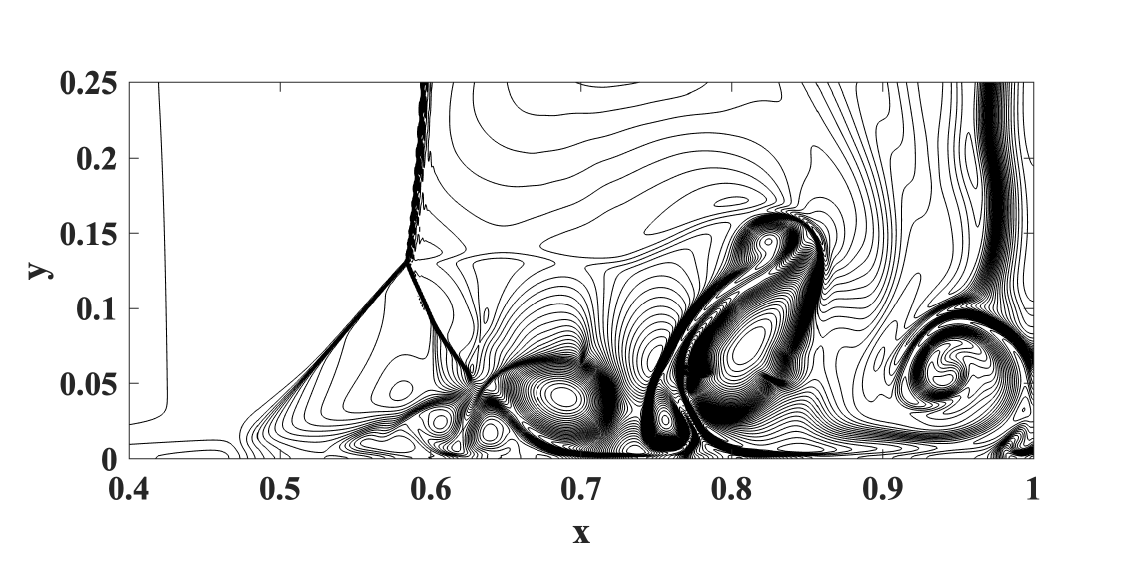}
	}
	\subfigure[ASE-DF(9,7,5,3), $\Delta x=\Delta y=1/500$]{
		\includegraphics[width=0.42\textwidth]{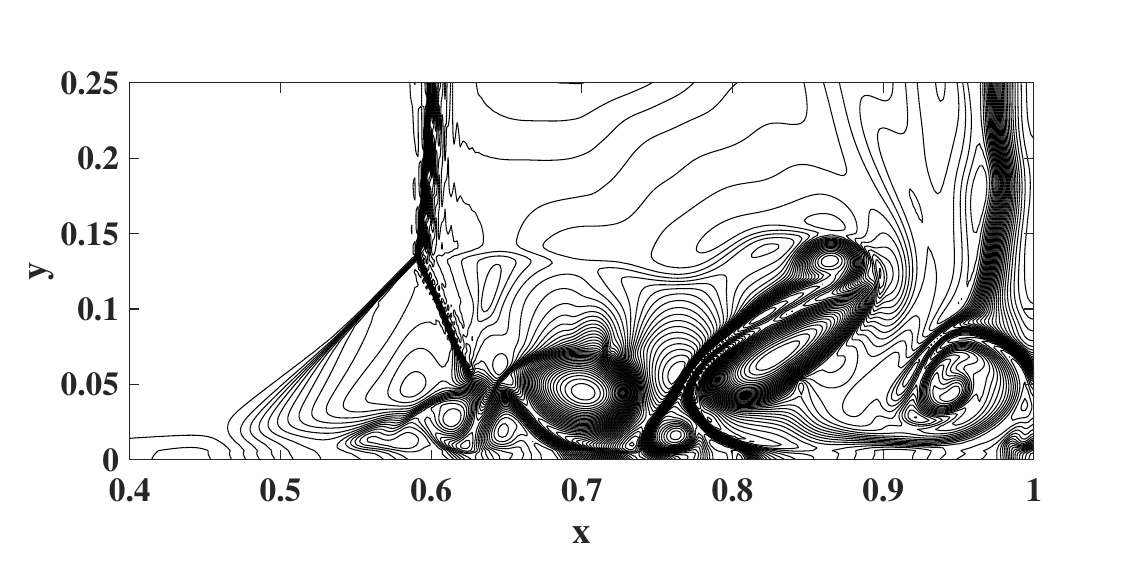}
	}
	\subfigure[ASE-DF(9,7,5,3), $\Delta x=\Delta y=1/1000$]{
		\includegraphics[width=0.42\textwidth]{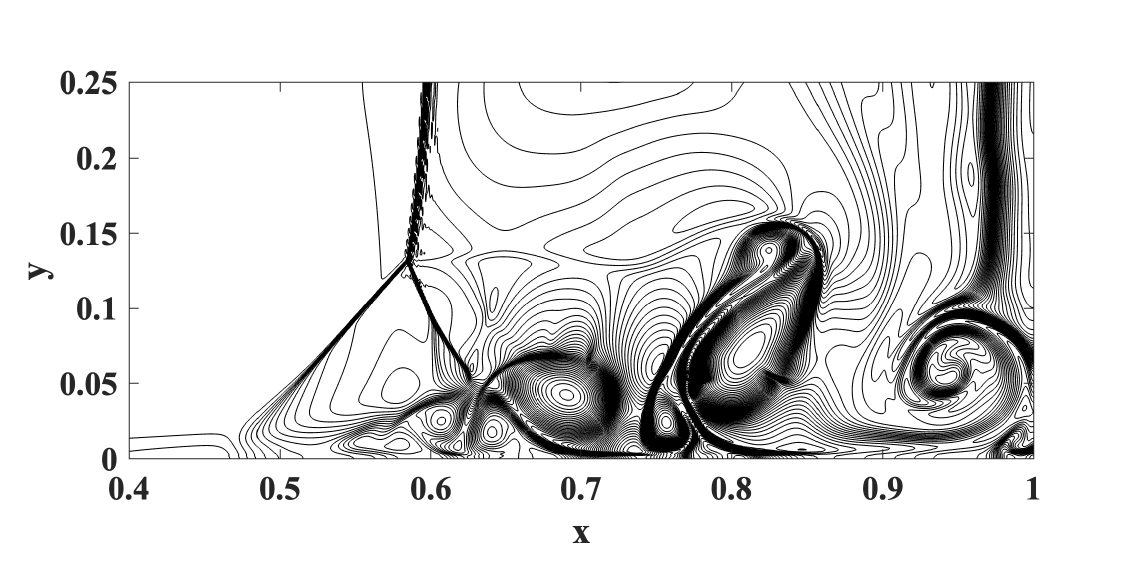}
	}
	\caption{\color{black}Viscous sod shock problem: the density distribution at $t=1.0$ for $Re=200$ case. Left: different reconstruct methods with 500$\times 250$ uniform meshes. Right: different reconstruct methods with 1000$\times 500$ uniform meshes. All results are calculated by the GKS solver, $c_1=1.0, c_2=10.0$.\color{black}}
	\label{viscous sod problem}
\end{figure}
\begin{figure}[htbp]
	\centering
	\subfigure[ASE-DF(5,3)]{
		\includegraphics[width=0.42\textwidth]{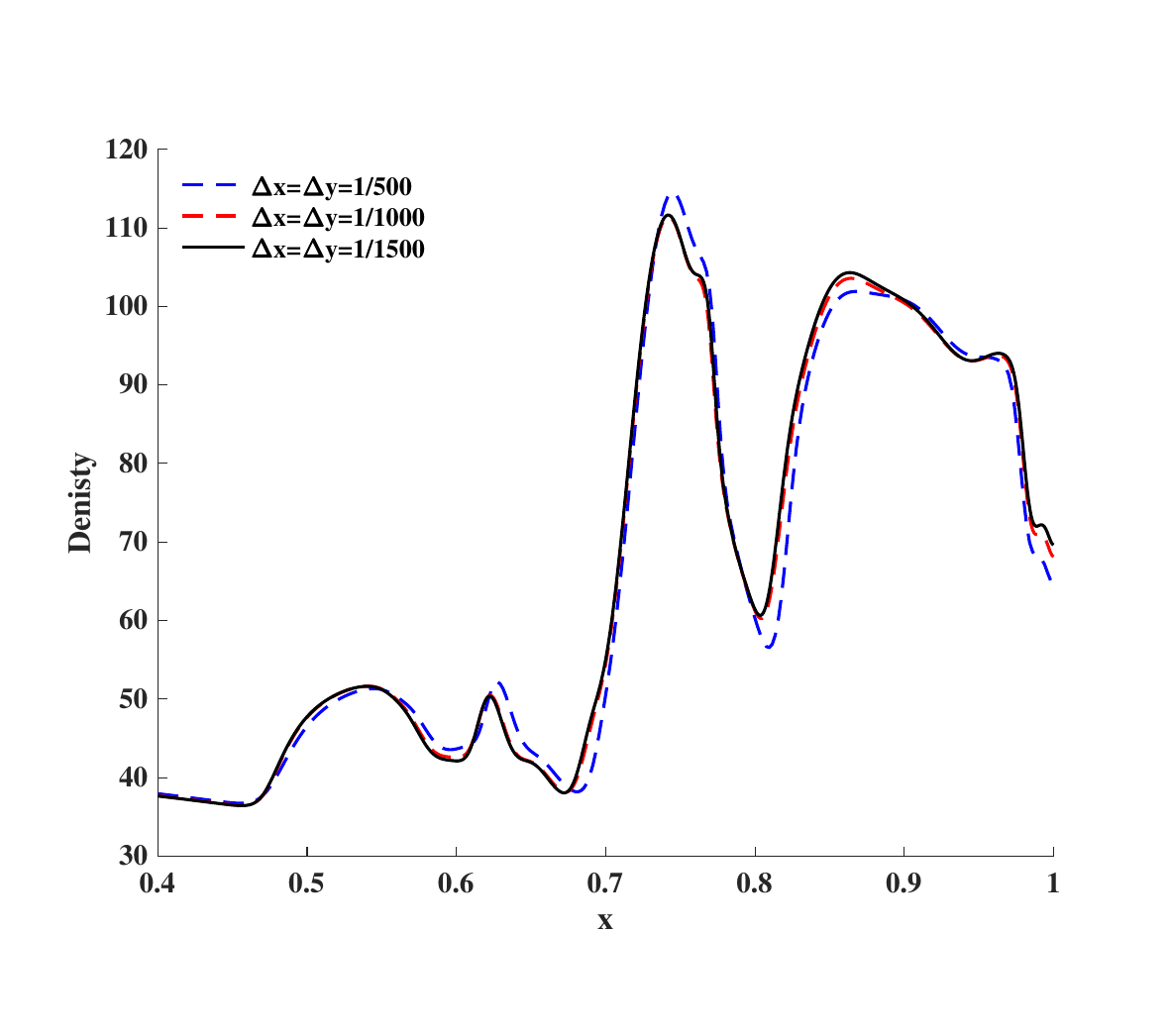}
	}
	\subfigure[ASE-DF(7,5,3)]{
		\includegraphics[width=0.42\textwidth]{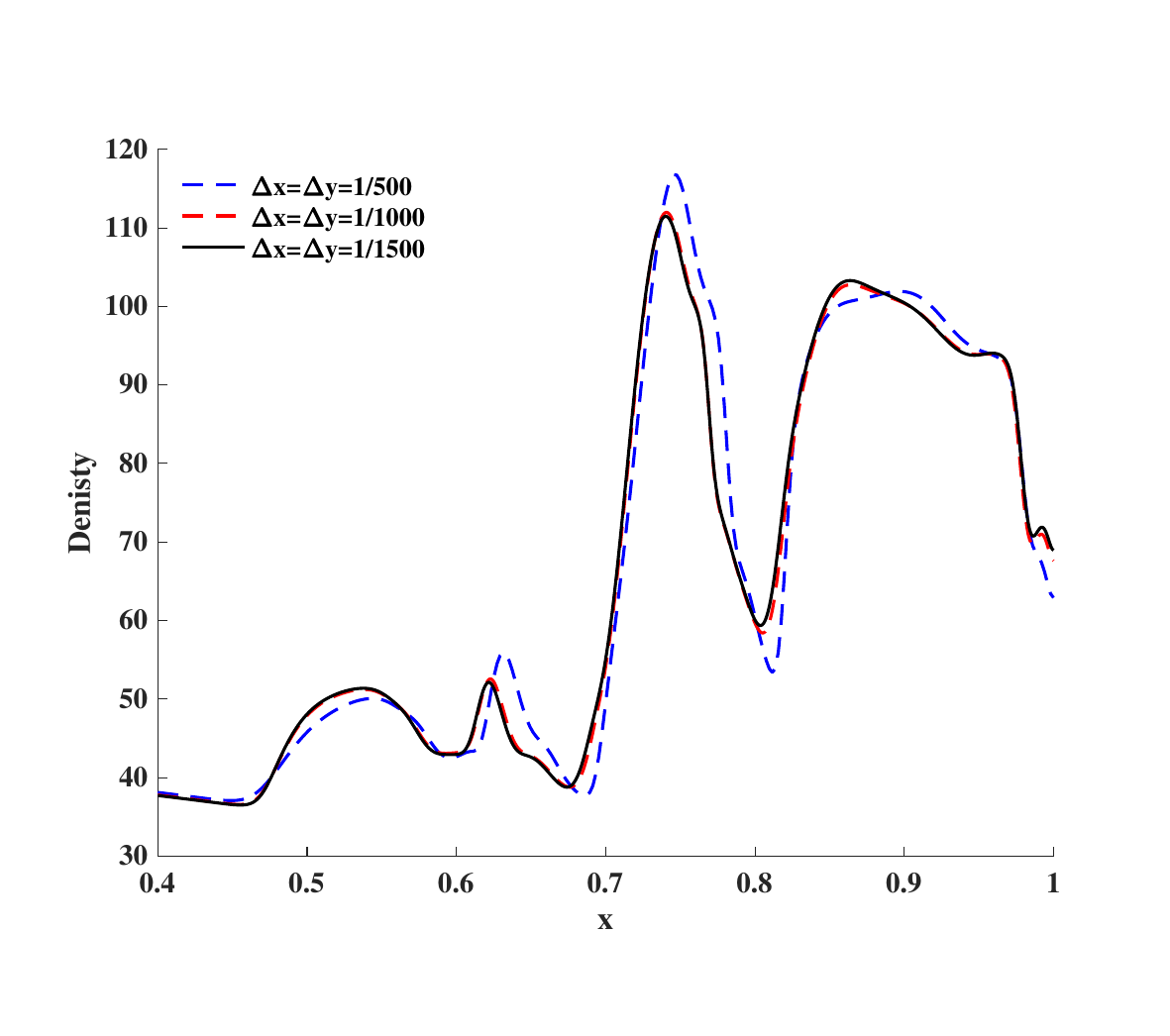}
	}
	\subfigure[ASE-DF(9,7,5,3)]{
		\includegraphics[width=0.42\textwidth]{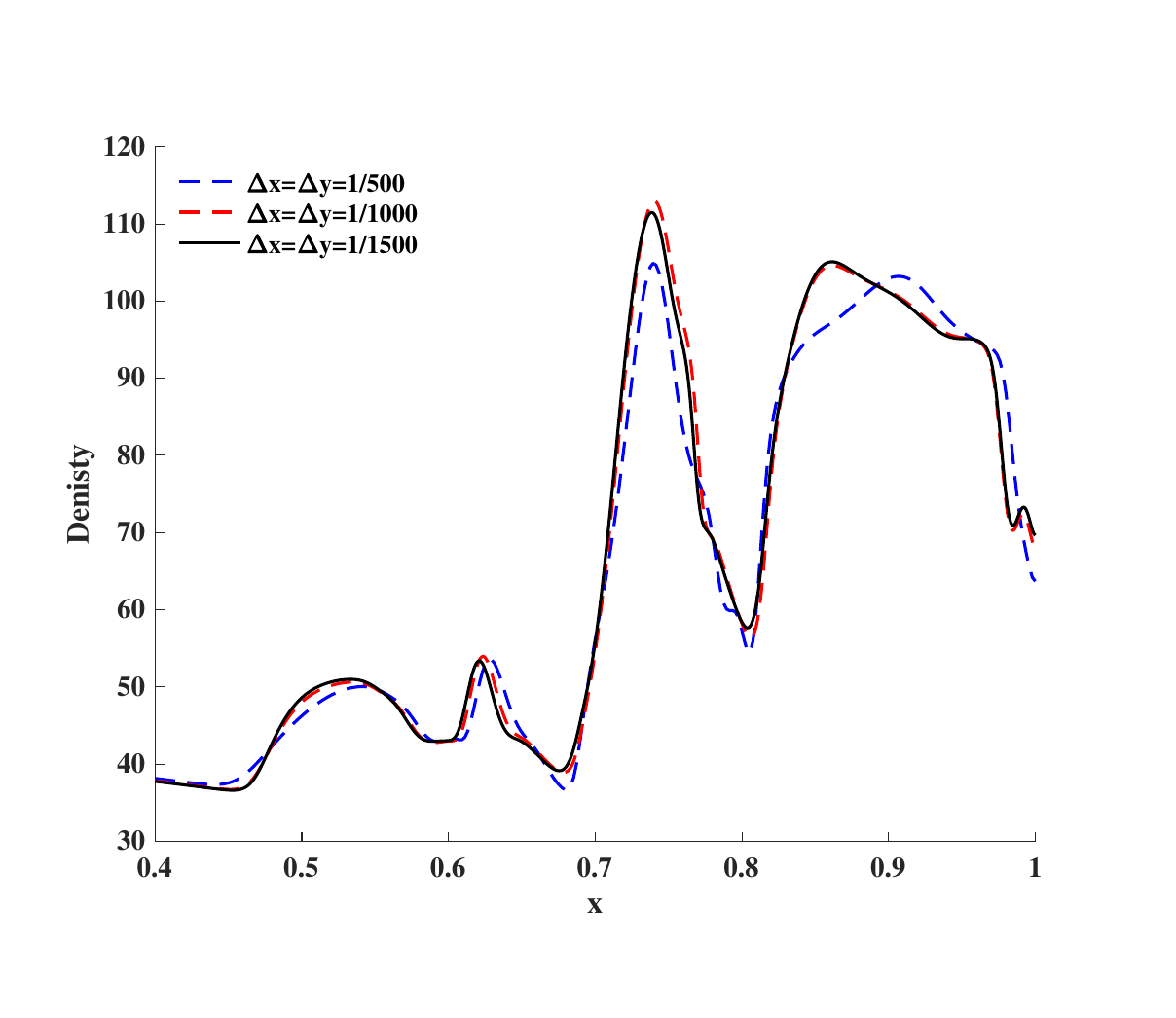}
	}
	\subfigure[Comparsion results]{
		\includegraphics[width=0.42\textwidth]{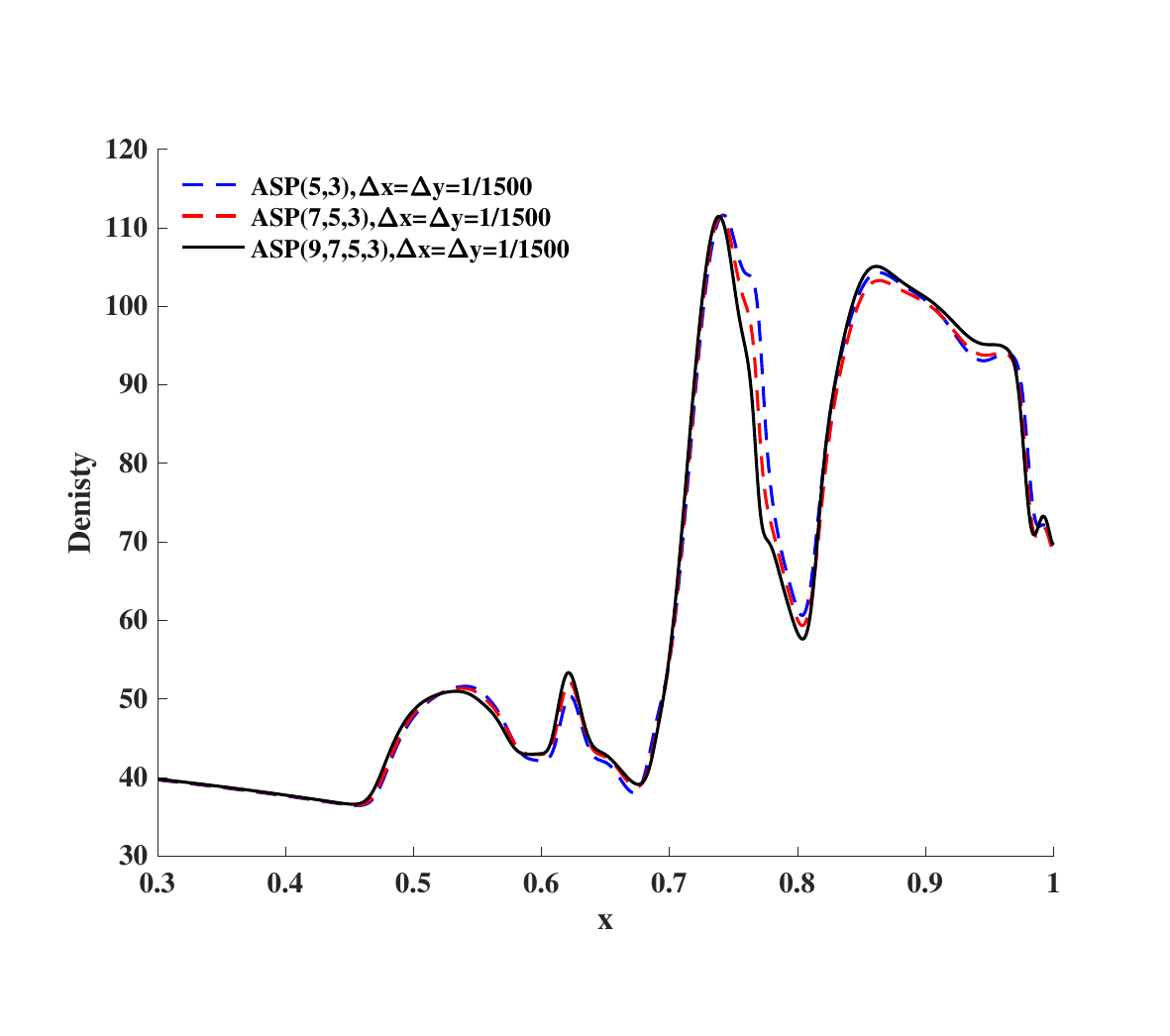}
	}
	\caption{Viscous sod shock problem: density profiles along the lower wall at $t=1.0$ for $Re=200$ case. (a-c) Mesh refinement results for different reconstruct methods. (d) Comparsion of different reconstruct methods with 1500$\times$750 uniform meshes.}
	\label{viscous sod problem-re200}
\end{figure}

\paragraph{Example 5}(Viscous sod shock problem)
For N-S solver, the viscous shock tube problem \cite{daru2000evaluation} is considered to test the validity of the scheme. The initial condition in the domain $[0,1]\times[0,0.5]$ is given by
\begin{equation*}
	(\rho,u,v,p)=
	\begin{cases}
		(120,0,0,120/\gamma),\ 0<x<0.5,\\
		(1.2,0,0,1.2/\gamma),\ 0.5\le x<1,\\
	\end{cases}
\end{equation*}
where $\gamma=1.4$ and the Prandtl number $Pr=1$. The top boundary is set as a symmetric boundary condition and no-slip adiabatic condition is used for others.  For the case $Re=200$, the density distributions with different uniform mesh points at $t=1.0$ from different reconstruct methods are shown in Fig \ref{viscous sod problem}. The density profiles along the lower wall for this case are presented in Fig \ref{viscous sod problem-re200}. As shown in Tab \ref{comparsion with reference}, the height of primary vortex predicted by the current schemes agree well with the reference data \cite{2005Accurate}. 

\begin{figure}[htbp]
	\centering
	\subfigure[ASE-DF(5,3),\ Mach 80]{
		\includegraphics[width=0.44\textwidth]{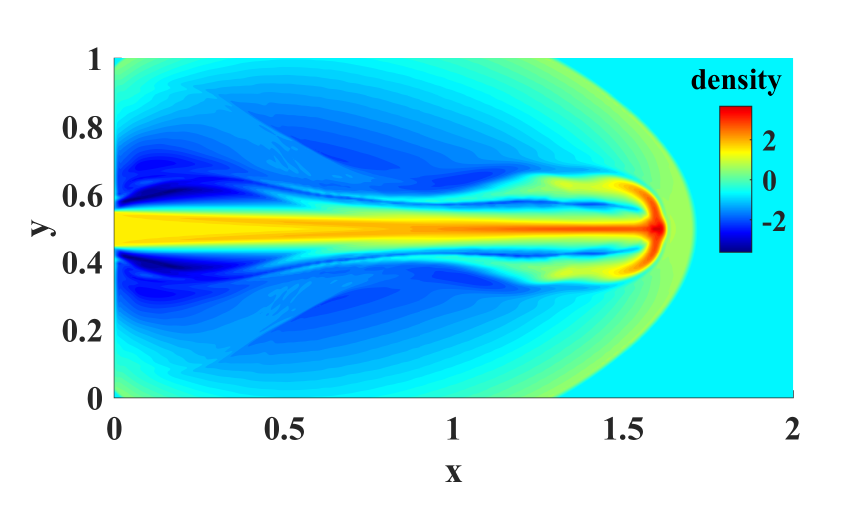}
	}
	\subfigure[ASE-DF(5,3),\ Mach 20000]{
		\includegraphics[width=0.3\textwidth]{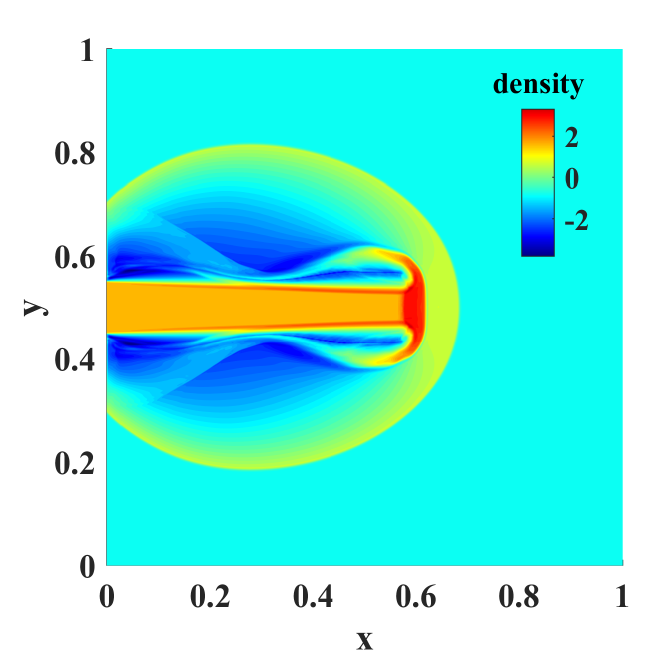}
	}\\
	\subfigure[ASE-DF(7,5,3),\ Mach 80]{
		\includegraphics[width=0.44\textwidth]{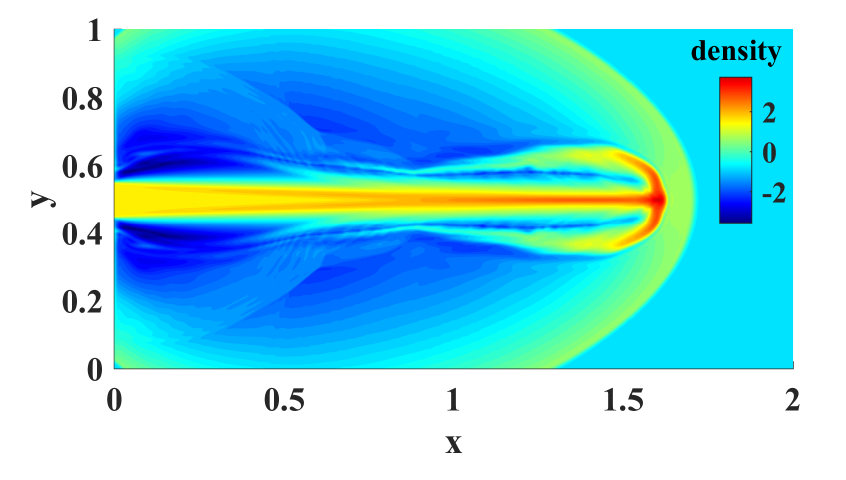}
	}
	\subfigure[ASE-DF(7,5,3),\ Mach 20000]{
		\includegraphics[width=0.3\textwidth]{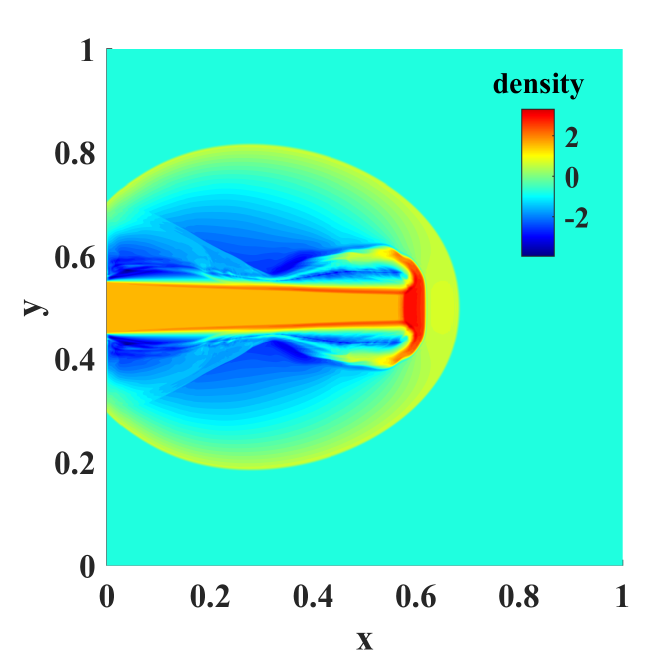}
	}\\
	\subfigure[ASE-DF(9,7,5,3),\ Mach 80]{
		\includegraphics[width=0.44\textwidth]{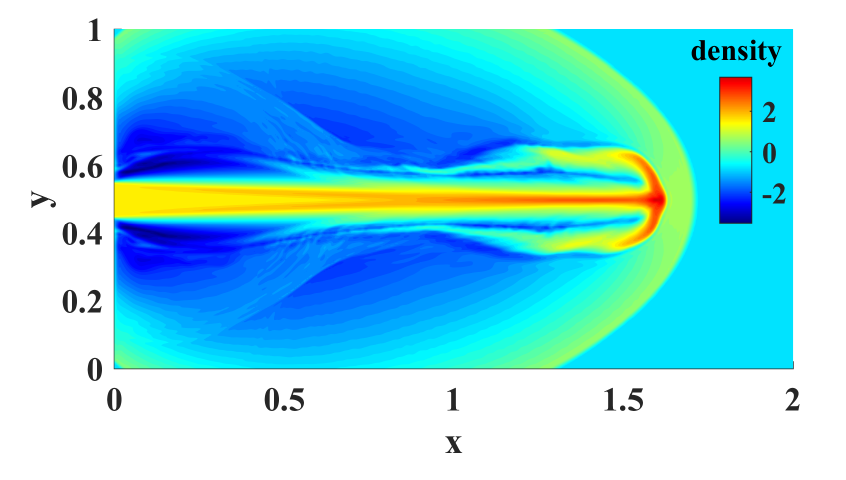}
	}
	\subfigure[ASE-DF(9,7,5,3),\ Mach 20000]{
		\includegraphics[width=0.3\textwidth]{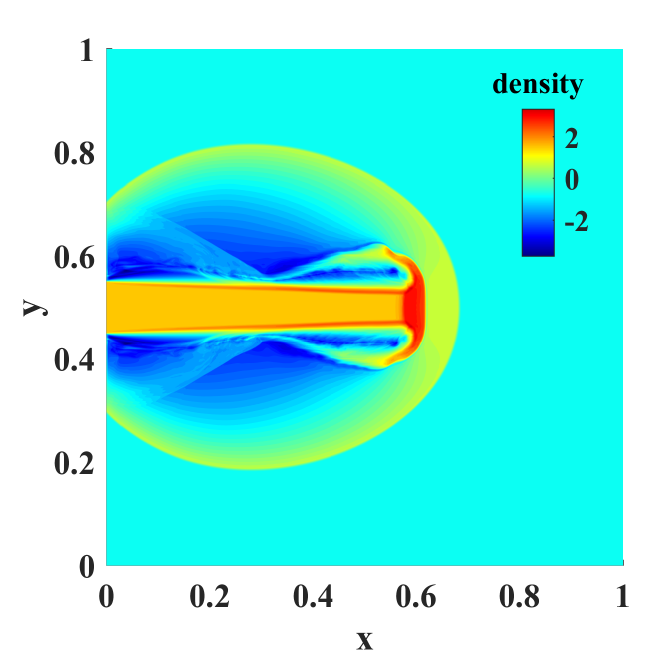}
	}
	\caption{High-mach number astrophysical jet: the density distribution at $t=0.07$ and $t=10^{-4}$(from left to right). This figure is drawn with 30 density contours. Non-linear function of density $\phi={\rm log}(\rho)$ is used. Left: different reconstruction methods with $400\times200$ meshes. Right: different reconstruction methods with $800\times400$ meshes. L-F solver is used for all results.}
	\label{high mach jet}
\end{figure}

\paragraph{Example 6}(High-mach number astrophysical jet) To test the robustness of the schemes, the high-mach number astrophysical jet \cite{Zhang2010On} is simulated. A low mach number $Ma\approx80$ and a very high mach number $Ma\approx20000$ are considered. The computational domain $[0,2]\times[0,1]$ is filled with
\begin{equation*}
	(\rho,u,v,p,\gamma)=(0.5,0,0,0.4127,5/3).	
\end{equation*}
The boundary conditions for the right, top and bottom are outflow. For the left boundary
\begin{equation*}
	\begin{aligned}
		(\rho,u,v,p)=(5,30,0,0.4127),\ {\rm if}\ y\in[0.45,0.55]\quad &{\rm For\ Ma\ 80\ case},\\
		(\rho,u,v,p)=(5,8000,0,0.4127),\ {\rm if}\ y\in[0.45,0.55]\quad &{\rm For\ Ma\ 20000\ case}.
	\end{aligned}
\end{equation*}
For the Ma 80 case, the terminal time is 0.07, the computation is performed on a $400\times200$ mesh, and for the Ma 20000 case, the terminal time is $10^{-4}$, the computation is performed on a $800\times400$ mesh. Fig \ref{high mach jet} shows the results when using different reconstruction methods, we see that the proposed algorithm can maintain high resolution at lower Mach number, while the algorithm maintains strong robustness at a very high Mach number.

\begin{table}[htbp]
	\caption{Viscous shock tube problem: comparison of the primary vortex heights among different schemes for Re = 200 case. The reference data can be found in \cite{2005Accurate}. }
	\label{comparsion with reference}       
	\color{black}
	\begin{tabular*}{\textwidth}{@{\extracolsep{\fill}}ccccc}
		\hline\noalign{\smallskip}
		Scheme &AUSMPW+&M-AUSMPW+ & ASP-DF(5,3)   & ASP-DF(7,5,3) \\
		\noalign{\smallskip}\hline\noalign{\smallskip}
		Mesh size &500$\times$250 &500$\times$250 &500$\times$250  &500$\times$250    \\
		Height &0.163 &0.168 &0.167  &0.163   \\
		\hline\noalign{\smallskip}
		Scheme &ASP-DF(9,7,5,3)& ASP-DF(5,3)   & ASP-DF(7,5,3)  & ASP-DF(9,7,5,3) \\
		\noalign{\smallskip}\hline\noalign{\smallskip}
		Mesh size &500$\times$250 &1000$\times$500 &1000$\times$500  &1000$\times$500    \\
		Height &0.153 &0.171 &0.164  &0.160  \\
		\noalign{\smallskip}\hline
	\end{tabular*}
\end{table}

\subsection{Computational efficiency}
The 2D Riemann problem (Configuration 3 in \cite{Lax1998Solution}) with different mesh sizes are tested for the comparison of computation cost, which is shown in Tab \ref{time cost}. The CPU times are recorded after running 20 time steps for each reconstruction method with a single processor of AMD EPYC 7401 @2.00GHz.

Based on the Tab \ref{time cost}, since there is no need to calculate additional smoothness indicators in higher-order schemes, \color{black}the computational time of ASE-DF(7,5,3) is about 1.4 times of ASE-DF(5,3) due to the differences in number of Gaussian points and the reconstruct polynomials, and the computational time of ASE-DF(9,7,5,3) is about 1.8 times of ASE-DF(5,3), \color{black} i.e. our proposed algorithm remains efficient at arbitrary high-order methods.

\begin{table}[h]
	\caption{Computational time (in seconds) of different reconstruction methods with the GKS solver for the 2-D Riemann problem.}
	\label{time cost}       
	\begin{tabular*}{\textwidth}{@{\extracolsep{\fill}}cccc}
		\hline\noalign{\smallskip}
		Mesh size & ASP-DF(5,3)   & ASP-DF(7,5,3)  & ASP-DF(9,7,5,3) \\
		\noalign{\smallskip}\hline\noalign{\smallskip}
		100$\times$100 &5.974  &9.053 &10.803   \\
		200$\times$200 &23.116  &32.852  &42.554   \\
		300$\times$300 &51.349  &72.788  &94.64  \\
		400$\times$400 &90.954  &128.006  &165.221 \\
		\noalign{\smallskip}\hline
	\end{tabular*}
\end{table}

\section{Conclusion}
This paper introduces an efficient class of adaptive stencil extension reconstruction methods based on a Discontinuity Feedback factor (ASE-DF). 
These methods can be readily extended to arbitrary high-order forms, addressing and releasing two key challenges faced by high-order WENO schemes:
Declining robustness in extreme problem calculations as order increases; Escalating computational costs due to additional smoothness indicator calculations.
Our approach builds upon the WENOZ-AO(5,3) scheme, incorporating the DF factor to significantly enhance algorithm robustness. Notably, this method eliminates the need for additional smoothness indicators, relying solely on the stencil's $\alpha_\mathbb{S}$ for expansion to higher orders. This ensures the algorithm maintains high resolution and efficiency.
The ASE-DF method offers a promising foundation for high-order schemes, preserving accuracy while mitigating computational overhead. Future research will explore the potential of this technique in compact schemes and multi-component flows.

\section*{Acknowledgments}
The current research is supported by the National Science Foundation of China (12302378,
12172316, 92371201), Hong Kong Research Grant Council (16208021).
\begin{appendix}
	\section{Calculation of GKS flux function}\label{GKS flux calculation}
	The form of the 2nd-order gas kinetic distribution function at each Gaussian point with a local coordinate $\mathbf{x}$ in Eq~(\ref{GKS flux})
	For a 2-D Maxwell's distribution
	\begin{equation}
		g=\rho\left(\frac{\lambda}{\pi}\right)^{\frac{K+3}{2}}e^{-\lambda\left[(u-U)^2+(v-V)^2+\xi^2\right]},
	\end{equation}
	the moment of $g$ is defined as
	\begin{equation}
		\rho\left<(\cdots)\right>=\int(\cdots)g{\rm d}\Xi,
		\label{moment function}
	\end{equation}
	the general moment formula becomes
	\begin{equation*}
		\left<(u^nv^m\xi^{2l})\right>=\left<(u^n)\right>\left<(v^m)\right>\left<(\xi^{2l})\right>,
	\end{equation*}
	where $n,m,l$ are integers, and the moments of $\xi$ are always even-order because of its symmetrical property. With the integral in the domain $(-\infty,+\infty)$, we have
	\begin{equation*}
		\begin{aligned}
			\left<(u^0)\right>&=1,\\
			\left<(u^1)\right>&=U,\\
			&\cdots\\
			\left<(u^{n+2})\right>&=U\left<(u^{n+1})\right>+\frac{n+1}{2\lambda}\left<(u^{n})\right>,\\
		\end{aligned}
	\end{equation*}
	and
	\begin{equation*}
		\begin{aligned}
			\left<(\xi^0)\right>&=1,\\
			\left<(u^2)\right>&=\frac{K}{2\lambda},\\
			&\cdots\\
			\left<(\xi^{2l})\right>&=\frac{K+2(l-1)}{2\lambda}\left<(\xi^{2(l-1)})\right>.\\
		\end{aligned}
	\end{equation*}
	Integrating terms with Heaviside function, the integral in the domain $(0,+\infty)$ is denoted as $\left<(\cdots)\right>_{>0}$, and $(-\infty,0)$ as $\left<(\cdots)\right>_{<0}$
	\begin{equation*}
		\begin{aligned}
			\left<(u^0)\right>_{>0}&=\frac12{\rm erfc}\left(-\sqrt{\lambda}U\right),\\
			\left<(u^1)\right>_{>0}&=U\left<(u^0)\right>_{>0}+\frac12\frac{e^{-\lambda U^2}}{\sqrt{\pi \lambda}},\\
			&\cdots\\
			\left<(u^{n+2})\right>_{>0}&=U\left<(u^{n+1})\right>_{>0}+\frac{n+1}{2\lambda}\left<(u^{n})\right>_{>0},\\
		\end{aligned}
	\end{equation*}
	and
	\begin{equation*}
		\begin{aligned}
			\left<(u^0)\right>_{<0}&=\frac12{\rm erfc}\left(\sqrt{\lambda}U\right),\\
			\left<(u^1)\right>_{<0}&=U\left<(u^0)\right>_{<0}-\frac12\frac{e^{-\lambda U^2}}{\sqrt{\pi \lambda}},\\
			&\cdots\\
			\left<(u^{n+2})\right>_{<0}&=U\left<(u^{n+1})\right>_{<0}+\frac{n+1}{2\lambda}\left<(u^{n})\right>_{<0},\\
		\end{aligned}
	\end{equation*}
	where erfc is the standard complementary error function.
	
	For the Taylor expansion of a Maxwell's distribution, all microscopic derivatives have the form
	\begin{equation*}
		\begin{aligned}
			&a_x=a_{x1}+a_{x2}u+a_{x3}v+a_{x4}\frac12(u^2+v^2+\xi^2)=a_{xi}\psi_i,\\
			&a_y=a_{y1}+a_{y2}u+a_{y3}v+a_{y4}\frac12(u^2+v^2+\xi^2)=a_{yi}\psi_i,\\
			&A=A_{1}+A_{2}u+A_{3}v+A_{4}\frac12(u^2+v^2+\xi^2)=A_{i}\psi_i.\\
		\end{aligned}
	\end{equation*}
	Combine with the Eq~(\ref{moment function}) and Eq~(\ref{relationship between g and W}), we have
	\begin{equation}
		\int \bm{\psi}a_xg{\rm d}\Xi=\frac{\partial\mathbf{W}}{\partial x},
	\end{equation}
	which can be expanded as
	\begin{equation}
		\left(\begin{array}{c}
			b_1  \\
			b_2\\
			b_3\\
			b_4
		\end{array}\right)=\frac1{\rho}\frac{\partial \mathbf{W}}{\partial x}=\frac1{\rho}
		\left(\begin{array}{c}
			\frac{\partial\rho}{\partial x}  \\
			\frac{\partial(\rho U)}{\partial x}\\
			\frac{\partial(\rho V)}{\partial x}\\
			\frac{\partial(\rho E)}{\partial x}
		\end{array}\right)=\left<(a_{xi}\psi_i\psi_j)\right>=\left<(\psi_i\psi_j)\right>\left(\begin{array}{c}
			a_{x1}  \\
			a_{x2}\\
			a_{x3}\\
			a_{x4}
		\end{array}\right).
	\end{equation}
	Denoting $\mathbf{M}=\left<(\psi_i\psi_j)\right>$, the above equations can be expressed as
	\begin{equation*}
		\mathbf{Ma=b}.
	\end{equation*}
	For the 2-D case, the coefficient matrix $\mathbf{M}$ can be expanded as
	\begin{equation}
		\mathbf{M}=
		\left(\begin{array}{cccc}
			\left<(u^0)\right>&\left<(u^1)\right>&\left<(v^1)\right>&\left<(\psi_4)\right>  \\
			\left<(u^1)\right>&\left<(u^2)\right>&\left<(u^1v^1)\right>&\left<(u^1\psi_4)\right>  \\
			\left<(v^1)\right>&\left<(u^1v^1)\right>&\left<(v^2)\right>&\left<(v^1\psi_4)\right>  \\
			\left<(\psi_4)\right>&\left<(u^1\psi_4)\right>&\left<(v^1\psi_4)\right>&\left<(\psi_4^2)\right>  \\
		\end{array}\right)=
		\left(\begin{array}{cccc}
			1 & U & V & B_1  \\
			U & U^2+1/2\lambda & UV & B_2  \\
			V & UV & V^2+1/2\lambda & B_3  \\
			B_1 & B_2 & B_3 & B_4  \\
		\end{array}\right),
		\label{coefficient}
	\end{equation}
	where
	\begin{equation*}
		\begin{aligned}
			B_1&=\frac12\left(U^2+V^2+\frac{K+2}{2\lambda}\right),\\
			B_2&=\frac12\left(U^3+V^2U+\frac{(K+4)U}{2\lambda}\right),\\
			B_3&=\frac12\left(V^3 + U^2V+\frac{(K+4)V}{2\lambda}\right),\\
			B_4&=\frac14\left((U^2+V^2)^2+\frac{(K+4)(U^2+V^2)}{\lambda}+\frac{K^2+6K+8}{4\lambda^2}\right).
		\end{aligned}
	\end{equation*}
	Denoting
	\begin{equation*}
		R_4=2b_4-\left(U^2+V^2+\frac{K+2}{2\lambda}\right)b_1,\quad R_3=b_3-Vb_1,\quad R_2=b_2-Ub_1,
	\end{equation*}
	the solution of Eq~(\ref{coefficient}) can be expressed as
	\begin{equation*}
		\begin{aligned}
			a_{x4}&=\frac{4\lambda^2}{K+2}(R_4-2UR_2-2VR_3),\\
			a_{x3}&=2\lambda R_3-Va_{x4},\quad a_{x2}=2\lambda R_2-U a_{x4},\\
			a_{x1}&=b_1-Ua_{x2}-Va_{x3}-\frac12a_{x4}\left(U^2+V^2+\frac{K+2}{2\lambda}\right).
		\end{aligned}
	\end{equation*}
	
\end{appendix}

\end{document}